\documentstyle[12pt]{amsart}

\numberwithin{equation}{section}

\begin{document}

%



\setlength{\itemsep}{0in}\newcommand{\lab}{\label}
\newcommand{\labeq}[1]{  \be \label{#1}  }
\newcommand{\labea}[1]{  \bea \label{#1}  }
\newcommand{\ben}{\begin{enumerate}}
\newcommand{\een}{\end{enumerate}}
\newcommand{\bm}{\boldmath}
\newcommand{\Bm}{\Boldmath}
\newcommand{\bea}{\begin{eqnarray}}
\newcommand{\ba}{\begin{array}}
\newcommand{\bean}{\begin{eqnarray*}}
\newcommand{\ea}{\end{array}}
\newcommand{\eea}{\end{eqnarray}}
\newcommand{\eean}{\end{eqnarray*}}
\newcommand{\beq}{\begin{equation}}
\newcommand{\eeq}{\end{equation}}
\newcommand{\bthm}{\begin{thm}}
\newcommand{\ethm}{\end{thm}}
\newcommand{\blem}{\begin{lem}}
\newcommand{\elem}{\end{lem}}
\newcommand{\bprop}{\begin{prop}}
\newcommand{\eprop}{\end{prop}}
\newcommand{\bcor}{\begin{cor}}
\newcommand{\ecor}{\end{cor}}
\newcommand{\bdfn}{\begin{dfn}}
\newcommand{\edfn}{\end{dfn}}
\newcommand{\brem}{\begin{rem}}
\newcommand{\erem}{\end{rem}}
\newcommand{\bex}{\begin{example}}
\newcommand{\eex}{\end{example}}
\newcommand{\lb}{\linebreak}
\newcommand{\nlb}{\nolinebreak}
\newcommand{\nl}{\newline}
\newcommand{\hs}{\hspace}
\newcommand{\vs}{\vspace}
\alph{enumii} \roman{enumiii}
\newtheorem{thm}{Theorem}[section]
\newtheorem{prop}[thm]{Proposition}
\newtheorem{lem}[thm]{Lemma}
\newtheorem{sublem}[thm]{Sublemma}
\newtheorem{cor}[thm]{Corollary}
\newtheorem{dfn}[thm]{Definition}
\newtheorem{rem}[thm]{Remark}
\newtheorem{defn}[thm]{Definition} 
\newtheorem{example}[thm]{Example} 

\def\endpf{\qed}
\def\ms{\medskip}
\def\N{{I\!\!N}}                \def\Z{Z\!\!\!\!Z}      \def\R{I\!\!R}
\def\C{{C\!\!\!\!I}}           \def\T{T\!\!\!\!T}      
\def\oc{\overline \C}
\def\Q{Q\!\!\!\!|}
\def\1{1\!\!1}

\def\and{\text{ and }}        \def\for{\text{ for }}
\def\tif{\text{ if }}         \def\then{\text{ then }}

\def\Cap{\text{Cap}}          \def\Con{\text{Con}}
\def\Comp{\text{Comp}}        \def\diam{\text{\rm {diam}}}
\def\dist{\text{{\rm dist}}}  \def\Dist{\text{{\rm Dist}}}
\def\Crit{\text{Crit}}
\def\Sing{\text{Sing}}        \def\conv{\text{{\rm conv}}}
\def\CIFS{\text{\rm CIFS}}    \def\PIFS{\text{\rm PIFS}}
\def\SIFS{\text{\rm SIFS}}
\def\IR{\text{\rm IR}}        \def\RE{\text{\rm R}}
\def\CR{\text{\rm CR}}        \def\SR{\text{\rm SR}}
\def\CFR{\text{\rm CFR}}

\def\Log{\text{\rm Log}}

\def\supp{\text{{\rm supp }}}

\def\Fin{{\cal F}in}
\def\F{{\Cal F}}
\def\h{{\text h}}
\def\hmu{\h_\mu}           \def\htop{{\text h_{\text{top}}}}

\def\H{\text{{\rm H}}}     \def\HD{\text{{\rm HD}}}   
\def\DD{\text{DD}}
\def\BD{\text{{\rm BD}}}         \def\PD{\text{PD}}
\def\re{\text{{\rm Re}}}    \def\im{\text{{\rm Im}}}
\def\Int{\text{{\rm Int}}} \def\ep{\text{e}}
\def\CD{\text{CD}}         \def\P{\text{{\rm P}}}     
\def\Id{\text{{\rm Id}}}
\def\Hyp{\text{{\rm Hyp}}}
\def\A{\Cal A}             \def\Ba{\Cal B}       \def\Ca{\Cal C}
\def\Ha{{\cal H}}
\def\LC{{\cal L}}             \def\M{\Cal M}        \def\Pa{\Cal P}
\def\U{\Cal U}             \def\V{\Cal V}
\def\W{\Cal W}

\def\a{\alpha}                \def\b{\beta}             \def\d{\delta}
\def\De{\Delta}               \def\e{\varepsilon}          
\def\g{\gamma}                \def\Ga{\Gamma}           
\def\L{\Lambda}              \def\l{\lambda}           
\def\Om{\Omega}               \def\om{\omega}
\def\Sg{\Sigma}               \def\sg{\sigma}
\def\Th{\Theta}               \def\th{\theta}           
\def\vth{\vartheta}
\def\Ka{\Kappa}               \def\ka{\kappa}

\def\bi{\bigcap}              \def\bu{\bigcup}
\def\({\bigl(}                \def\){\bigr)}
\def\lt{\left}                \def\rt{\right}
\def\bv{\bigvee}
\def\ld{\ldots}               \def\bd{\partial}         \def\^{\tilde}
\def\club{\hfill{$\clubsuit$}}\def\proot{\root p\of}

\def\tm{\widetilde{\mu}}
\def\tn{\widetilde{\nu}}
\def\es{\emptyset}            \def\sms{\setminus}
\def\sbt{\subset}             \def\spt{\supset}

\def\gek{\succeq}             \def\lek{\preceq}
\def\eqv{\Leftrightarrow}     \def\llr{\Longleftrightarrow}
\def\lr{\Longrightarrow}      \def\imp{\Rightarrow}
\def\comp{\asymp}
\def\upto{\nearrow}           \def\downto{\searrow}
\def\sp{\medskip}             \def\fr{\noindent}        
\def\nl{\newline}

\def\ov{\overline}            \def\un{\underline}

\def\ess{{\rm ess}}
\def\ni{\noindent}
\def\cl{\text{cl}}
\def\bt{{\bf t}}
\def\Bu{{\bf u}}
\def\tr{t}
\def\bo{{\bf 0}}
\def\nut{\nu_\lla^{\scriptscriptstyle 1/2}}
\def\arg{\text{arg}}
\def\Arg{\text{Arg}}
\def\re{\text{{\rm Re}}}
\def\gr{\nabla}
\def\endpf{{${\mathbin{\hbox{\vrule height 6.8pt depth0pt width6.8pt  
}}}$}}
\def\Fa{\cal F}
\def\Gal{\cal G}
\def\1{1\!\!1}
\def\D{{I\!\!D}}

\title{Multifractal analysis for conformal graph directed Markov systems}
\author[\sc Mario ROY]{\sc Mario ROY}
\address{Mario Roy, Glendon College, York University, 
2275 Bayview Avenue, Toronto, Canada, M4N 3M6;
\newline {\tt mroy@@gl.yorku.ca}}
\author[\sc Mariusz URBA\'NSKI]{\sc Mariusz URBA\'NSKI}
\address{Mariusz Urba\'nski, Department of Mathematics,
University of North Texas, Denton, TX 76203-1430, USA;
\newline {\tt urbanski@@unt.edu
\newline Webpage: www.math.unt.edu/$\sim$urbanski}}
%
\thanks{Research of the first author was supported by NSERC 
(Natural Sciences and Engineering Research Council of Canada).
Research of the second author was supported in part by the
NSF Grant DMS 0700831. Part of the work has been done while 
the second author was visiting the Max Planck Institute in Bonn, 
Germany. He wishes to thank the institute for its support.}

\begin{abstract}
We derive the multifractal analysis of the conformal measure 
(or equivalently, the invariant measure)
associated to a family of weights imposed upon a (multi-dimensional) graph 
directed Markov system 
(GDMS) using balls as the filtration. This analysis is done over 
a subset of $J$ which is often large. In particular, it coincides with
the limit set when the GDMS under scrutiny satisfies a boundary separation condition.
It also applies to more general situations such as real or complex continued fractions.
\end{abstract}

\maketitle

\section{Introduction}\label{intro}
Multifractal formalism origins from physics and mathematics 
(among others, see~\cite{mand}, \cite{fp}, \cite{gp} and~\cite{hjk}). 
In this latter paper, strong hints of parallels between multifractal
theory and the theory of statistical physics were suggested. Some of the 
first rigorous mathematical results on multifractals can be found in~\cite{cm} 
and~\cite{rand}. Since then, many papers have been written on this subject (for
instance, see~\cite{Ol1},~\cite{Ol2},~\cite{Pat} and~\cite{PW}). In particular, 
Pesin~\cite{pesin} developed a general framework in which multifractal formalism 
can be derived. 

We now briefly describe our setting. Let $\mu$ be a Borel probability measure on 
a metric space $X$. The measure $\mu$ is said to have local dimension $\a$ at a
point $x\in X$ if 
\[ \lim_{r\rightarrow 0}\frac{\log\mu(B(x,r))}{\log r}=\a. \] 
For each number $\a\geq 0$, let $X_\mu(\a)$ be the set of points $x\in X$ 
where the measure $\mu$ has local dimension $\a$, and let
$f_\mu(\a)$ be the Hausdorff dimension of the set $X_\mu(\a)$. 
The map $\a\mapsto f_\mu(\a)$ is called the (fine Hausdorff) multifractal 
spectrum of the measure $\mu$. 

The multifractal analysis of equilibrium 
states for a natural potential function and a
natural family of H\"older continuous weights
was performed in~\cite{HMU} for infinite conformal iterated function systems 
and in~\cite{gdms} for infinite conformal graph directed Markov systems. 
In both cases, the authors used at every point in the limit set of the given 
system the natural filtration generated by the initial blocks of the word that 
encodes the point. In other terms, the analysis was carried out using cylinders. 

Aiming to give the multifractal analysis a transparent geometrical meaning,
we shall derive in the sequel the multifractal analysis for cofinitely regular 
graph directed Markov systems (GDMSs) using as the filtration a base of balls 
centred at the given point. Until now, the question of the analysis of balls, 
which has already been solved in the case of finite systems, remained open 
for infinite systems. 

We conduct our analysis over the full set of parameters on which it 
can be expected to hold. Moreover, we perform it on a large and dynamically 
significant subset $J_r$ of the limit set $J$. In fact, under a mild boundary 
separation condition (for example, the separation condition) the sets $J_r$ and  
$J$ coincide. An additional geometric flavor of our analysis results from 
the fact that we concentrate on a geometrically meaningful family of 
H\"older continuous weights. All our results apply
to a large class of GDMSs, one- and multi-dimensional alike, including real and 
complex continued fractions.

Let us describe the content of the paper more precisely. In section~\ref{prelgdms}, 
we recall the basic definition of conformal graph directed Markov system. In 
section~\ref{hpt}, we describe the H\"older families of weights
$F$ and $F_{q,t}$ we shall work with and study the properties of 
the pressure $P(q,t)$ and temperature $T(q)$ they determine. In 
section~\ref{jr}, we carry out the multifractal analysis of the 
conformal measure $m_F$ (or equivalently, the invariant measure $\mu_F$) 
associated to a family of weights $F$. This analysis is done
over a subset $J_r$ of $J$ and conducted by means of balls. In particular, 
we show that for each $\a$ there is an auxiliary measure that witnesses the 
Hausdorff dimension of the set $J_{r,\mu}(\a)$ and that the $f_{r,\mu}(\a)$ 
curve is the Legendre transform of the temperature function $T(q)$. In 
section~\ref{jcond}, we derive the multifractal analysis of the conformal 
measure $m_F$ (or equivalently, the invariant measure $\mu_F$) 
under additional conditions on the GDMS. In subsection~\ref{ssc}, we 
observe that $J_r=J$ for all GDMSs which satisfy a boundary separation 
condition. Real continued fractions with the digit $1$ deleted are an 
important example of such systems. In subsection~\ref{oc}, we derive 
the multifractal analysis over $J$ under three conditions and show 
that there are families of one-dimensional conformal iterated function 
systems that meet these three conditions. Real continued fractions 
(with or without the digit $1$) are a good example of such a family. 
  
\section{Preliminaries on Graph Directed Markov Systems}\label{prelgdms}
Let us first describe the setting of conformal graph directed Markov systems 
introduced in~\cite{gdms}. Graph directed Markov systems are based on a 
directed multigraph $(V,E,i,t)$ and an associated incidence matrix $A$.
The multigraph consists of a finite set $V$ of vertices, a countable
(finite or infinite) set of edges, and two functions $i,t:E\rightarrow V$ that
indicate for each directed edge $e\in E$ its initial vertex $i(e)$ and its terminal vertex $t(e)$,
respectively. The matrix $A:E\times E\rightarrow\{0,1\}$ is an edge incidence matrix
and thus tells which edges may follow a given edge. Moreover, it respects the 
multigraph, that is, if $A_{e_1e_2}=1$ then $t(e_1)=i(e_2)$. It is thereafter natural to define
the set of all one-sided infinite $A$-admissible words
\[ E_A^\infty=\{\om\in E^\infty\,|\,A_{\om_i\om_{i+1}}=1, \forall\, i\in\N\}. \]
The set of all subwords of $E_A^\infty$ of length $n\in\N$ will be denoted by $E_A^n$, whereas 
the set of all finite subwords will be denoted by $E_A^*=\cup_{n\in\N}E_A^n$. 
From a dynamical point of view, we will consider the left shift map $\sg:E_A^\infty
\rightarrow E_A^\infty$ which drops the first letter of each word. 
 
A graph directed Markov system (GDMS) consists of a directed multigraph and an edge incidence matrix
together with a set of non-empty compact subsets $\{X_v\}_{v\in V}$ of a common Euclidean space $\R^d$, 
a number $0<s<1$, and for every $e\in E$ a one-to-one contraction $\varphi_e:X_{t(e)}\rightarrow X_{i(e)}$ 
with Lipschitz constant at most $s$. 

A GDMS is called iterated function system (IFS) provided that 
$A_{ef}=1$ if and only if $t(e)=i(f)$ and that $V$ is a singleton.

For $\om\in E_A^n$, $n\in\N$, we define
\[ \varphi_\om:=\varphi_{\om_1}\circ\varphi_{\om_2}\circ\cdots\circ
 \varphi_{\om_n}:X_{t(\om)}\rightarrow X_{i(\om)}. \] 
Note that the functions $i$ and $t$ extend naturally to $E_A^*$ by setting $i(\om):=i(\om_1)$
and $t(\om)=t(\om_{|\om|})$.
The main object of our interest will be the limit set $J$ of $S$. This set is the image of the symbolic space 
$E_A^\infty$ under a coding map $\pi$. Indeed, given any $\om\in E_A^\infty$, the sets $\varphi_{\om|_n}(X_{t(\om|_n)})$ 
form a decreasing sequence of non-empty compact sets whose diameters converge to zero. Therefore their intersection
$$
\bigcap_{n=1}^\infty\varphi_{\om|_n}(X_{t(\om|_n)})
$$
is a singleton, and we denote its element by $\pi(\om)$. This defines the coding map $\pi:E_A^\infty\to X$, 
where $X:=\bigoplus_{v\in V}X_v$ is the disjoint union of the compact sets $X_v$. Clearly, $\pi$ is a continuous 
function when $E_A^\infty$ is equipped with the topology generated by the cylinders 
$[e]_{n}=\{\om\in E_A^\infty\,|\,\om_n=e\}$, $e\in E$, $n\in\N$. 
Hence the limit set of the GDMS $S$ is  
$$
J=\pi(E^\infty)=\bigcup_{\om\in E^\infty}\bigcap_{n=1}^\infty\varphi_{\om|n}(X_{t(\om|_n)}).
$$

\ms\ni Recall also (cf. section~4.2 in~\cite{gdms}) that a GDMS $S=\{\varphi_e\}_{e\in E}$ is called conformal 
(and thereafter a CGDMS) if the following conditions are satisfied.
\begin{itemize}
\item[(i)] For every $v\in V$, the set $X_v$ is a compact, connected subset of $\R^d$ 
           which is the closure of its interior (i.e. $X_v=\overline{\Int_{\R^d}(X_v)}$);
\item[(ii)] (Open set condition (OSC)) For all $e,f\in E$, $e\neq f$,
\[ \varphi_e(\Int(X_{t(e)}))\cap\varphi_f(\Int(X_{t(f)}))=\emptyset; \]  
\item[(iii)] For every vertex $v\in V$, there exists an open connected set $W_v$ 
such that $X_v\sbt W_v\sbt\R^d$ and such that for every $e\in E$ with $t(e)=v$, the map 
$\varphi_e$ extends to a $C^1$ conformal diffeomorphism of $W_v$ into $W_{i(e)}$;
\item[(iv)] (Cone property) There exist $\g,l>0$ such that for every $v\in V$ and every $x\in X_v$ 
there is an open cone $\Con(x,\g,l)\sbt \Int(X_v)$ with vertex $x$,
central angle $\g$, and altitude $l$;
\item[(v)] There are two constants $L\ge 1$ and $\a>0$ such that
$$
\bigl||\varphi_e'(y)|-|\varphi_e'(x)|\bigr|\le
 L\|(\varphi_e')^{-1}\|^{-1}\|y-x\|^\a
$$
for every $e\in E$ and every pair of points $x,y\in W_{t(e)}$, where $|\varphi_e'(x)|$ 
denotes the norm of the derivative of $\varphi_e$ at $x$ and $\|(\varphi_e')^{-1}\|$ is the
supremum norm taken over $W_{i(e)}$.
\end{itemize}

\

\brem\lab{rsbdp}
According to Proposition~4.2.1 in \cite{gdms}, condition $(v)$ is 
automatically satisfied (with $\a=1$) when $d\ge 2$. This condition is also fulfilled if $d=1$, the alphabet 
$E$ is finite and all the $\varphi_e$'s are of class $C^{1+\e}$.
\erem

\

\ni The following useful fact has been proved in~Lemma~4.2.2 of~\cite{gdms}. 

\

\blem\lab{lsbdp}
For all $\om\in E_A^*$ and all $x,y\in W_{t(\om)}$ we have
$$
\bigl|\log|\varphi_\om'(y)|-\log|\varphi_\om'(x)|\bigr|
\le L(1-s)^{-1}\|y-x\|^\a.
$$
\elem

\

\ni An immediate consequence of this lemma is the famous bounded distortion property.

\

\begin{itemize}
\item[(v')] (Bounded Distortion Property (BDP)) There exists a constant
 $K\geq
1$ such that
$$
|\varphi_\om'(y)|\leq K|\varphi_\om'(x)|
$$
for every $\om\in E_A^*$ and every $x,y\in W_{t(\om)}$.
\end{itemize}

\

Recall that a CGDMS $S$ satisfy the Strong Open Set Condition (SOSC) if 
$J\cap\Int(X)\neq\emptyset$, that is, $\cup_{v\in V}(J_v\cap\Int(X_v))\neq\emptyset$.  
Recall further that a matrix $A$ is finitely primitive if there exists a finite set $\Om\sbt E_A^*$ of words 
of the same length such that for all $e,f\in E$ there is a word $\om\in\Om$ for which $e\om f\in E_A^*$. 
From this point on we assume that all the systems we deal with satisfy those two properties. 


\

Infinite systems naturally break into two main classes called 
irregular and regular systems. This
dichotomy can be determined from the existence of a conformal
measure or, equivalently, the existence of a zero of the 
topological pressure function. Recall that the topological pressure  
$P(t)$, $t\ge 0$, is 
defined as follows. For every $n\in\N$, set
$$
P^{(n)}(t)=\sum_{\om\in E_A^n}\|\varphi_\om'\|^t,
$$
where $\|\varphi_\om'\|:=\sup_{x\in X_{t(\om)}}|\varphi_\om'(x)|$.
Then
$$
P(t)=\lim_{n\to\infty}{1\over n}\log
P^{(n)}(t)=\inf_{n\in\N}{1\over n}\log P^{(n)}(t).
$$
If the function $\tilde{\zeta}:E_A^\infty\to\R$ is given by
the formula
$$
\tilde{\zeta}(\om)=\log|\varphi_{\om_1}'(\pi(\sg(\om)))|,
$$
then $P(t)=P(t\tilde{\zeta})$, where $\P(t\tilde{\zeta})$ is the
 classical topological pressure of the function $t\tilde{\zeta}$ when $E$ is finite 
(so the space $E^\infty$ is compact), and is understood 
in the sense of \cite{HMU} and \cite{gdms} when $E$ is infinite.
The finiteness parameter $\th$ of the system 
is defined by $\inf\{t\ge 0:P^{(1)}(t)<\infty\}=\inf\{t\ge 0:P(t)<\infty\}$. It 
is easy to show that the pressure function is non-increasing on $[0,\infty)$, 
that it is (strictly) decreasing, continuous and convex on 
$[\th,\infty)$, and that $P(d) \leq 0$. Of course,
 $P(0)=\infty$ if and 
only if $E$ is infinite.
The following characterization of the Hausdorff dimension 
$\HD(J)$ of the limit 
set $J$ was proved in \cite{gdms}, Theorem~4.2.13. 
For every $F\sbt E$, we write $S|_F$ for the 
subsystem $\{\varphi_e\}_{e\in F}$ of $S$, and $J_F$ for the limit set of $S|_F$.

\

\noindent{\bf Theorem.}
$$
\HD(J)=\inf\{t\ge 0:P(t)\le 0\}= \sup\{\HD(J_F):F\sbt E  \text{ is finite}\}\ge\th.
$$
If $P(t)=0$, then $t$ is the only zero of the function $P(t)$ and $t=\HD(J)$.

\

A system $S$ was called regular provided there is some 
$t\ge 0$ such that $P(t)=0$. In fact, a system is regular 
if and only if it admits a conformal measure. Recall that a Borel 
probability measure $m$ is said to be $t$-conformal provided $m(J)=1$
and for every $\om\in E_A^*$ and for every Borel set $B\sbt X_{t(\om)}$
$$
m(\varphi_\om(B))=\int_B|\varphi_\om'|^t\,dm, 
$$
and for all incomparable words $\om,\tau\in E_A^*$
\[ m\bigl(\varphi_\om(X_{t(\om)})\cap\varphi_\tau(X_{t(\tau)})\bigr)=0. \]

\sp\fr There are natural subclasses of regular systems. Among others,
%
%
%
%
%
a system is called
cofinitely regular provided every non-empty 
cofinite subsystem $S' = \{\varphi_e\}_{e\in E'}$ (i.e. $E'$ is a
cofinite subset of $E$) is regular.
A finite system is clearly cofinitely regular, and it 
was shown in \cite{gdms}, Theorem~4.3.4 that an infinite system is cofinitely regular
exactly when the pressure is infinite at the finiteness parameter.

\

\noindent{\bf Theorem.} 
An infinite system $S$ is cofinitely regular if and
only if $P(\th)=\infty \eqv P^{(1)}(\th)=\infty
\eqv \{t\ge 0:P(t)<\infty\}=(\th,\infty)
\eqv \{t\ge 0:P^{(1)}(t)<\infty\}=(\th,\infty)$.

%

\

Throughout the rest of this paper, all systems under investigation are assumed to be cofinitely regular.

\section{H\"older families of functions, Pressure and Temperature}\label{hpt}

We adopt the notational convention that families of functions shall be denoted by uppercase letters, 
while their members will be denoted by lowercase letters. Moreover, functions and measures associated 
with the symbolic space $E_A^\infty$ will wear a tilde $\sim$, with the notable exception of the shift 
map $\sg$.

Throughout this paper, let $S=\{\varphi_e:X_{t(e)}\rightarrow X_{i(e)}\,|\,e\in E\}$ be a cofinitely regular CGDMS satisfying the Strong Open Set Condition (SOSC) and having an underlying finitely primitive edge incidence matrix $A$. 
Let $\th$ be the finiteness parameter of $S$. Let $u>\th$. Let 
\[ 
\Psi=\bigl\{\psi_e:X_{t(e)}\rightarrow\R\,\bigl|\,e\in E\bigr\} 
\]
be a bounded H\"older family of functions. Let $\b$ be the order of that family. H\"older of order $\b$
means that (cf. section~3.1 in~\cite{gdms})
\[ 
V_\b(\Psi):=\sup_{n\in\N}V_n(\Psi)<\infty, 
\] 
where 
\[ 
V_n(\Psi)=\sup_{\om\in E_A^n}\sup_{x,y\in X_{t(\om)}}
          \bigl|\psi_{\om_1}(\varphi_{\sg\om}(x))-\psi_{\om_1}(\varphi_{\sg\om}(y))\bigr|e^{\b(n-1)}.
\]
Bounded simply means that 
\[ 
\|\Psi\|:=\sup_{e\in E}\|\psi_e\|<\infty. 
\]
Denote by $\tilde{\psi}:E_A^\infty\rightarrow\R$ the potential function (also called amalgamated function) 
induced by the family $\Psi$, which is defined by 
\[ 
\tilde{\psi}(\om)=\psi_{\om_1}(\pi(\sg\om)). 
\]  
According to Lemma~3.1.3 in~\cite{gdms}, the function $\tilde{\psi}$ is bounded and
H\"older continuous of order $\b$. Let also 
\[ 
\Log=\bigl\{\log|\varphi_e'|:X_{t(e)}\rightarrow\R\,\bigl|\,e\in E\bigr\}. 
\]
By Lemma~4.2.2 in~\cite{gdms}, the family $\Log$ is H\"older of order $\a\log s$. Moreover, for any $t>\th$, 
the family $t\Log$ is summable, that is
\[ 
\sum_{e\in E}\bigl\|\exp(t\log\|\varphi_e'\|)\bigr\|<\infty. 
\] 
In particular, $u\Log$ is summable.
Denote by $\tilde{\zeta}:E_A^\infty\rightarrow\R$ the amalgamated function induced by the family $\Log$, 
that is, 
\[ 
\tilde{\zeta}(\om)=\log|\varphi_{\om_1}'(\pi(\sg\om))|. 
\]
Then $\tilde{\zeta}$ is a summable H\"older continuous function of order $\a\log s$. It follows from the above 
definitions and properties that the family $F=\Psi+u\Log$, that is,
$F=\{f_e:X_{t(e)}\rightarrow\R\,|\,e\in E\}$,  
where 
\[ 
f_e=\psi_e+u\log|\varphi_e'|, 
\] 
is a summable H\"older family of functions of order 
$\g=\max\{\beta,\a\log s\}$. The amalgamated function $\tilde{f}:E_A^\infty\rightarrow\R$ 
induced by the family $F$ satisfies 
\[ 
\tilde{f}=\tilde{\psi}+u\,\tilde{\zeta} 
\]
and is a summable H\"older continuous function of order $\g$. 
Note also that $F$ and $\tilde{f}$ are bounded above by $\sup\Psi:=\sup_{e\in E}\sup_{x\in X_{t(e)}}\psi_e(x)$.
Recall that the topological pressure $P(F)$ of $F$ is defined by 
\[
P(F)=\lim_{n\rightarrow\infty}\frac{1}{n}\log\sum_{\om\in E_A^n}
         \exp\bigl(\sup_{x\in X_{t(\om)}}\sum_{i=1}^{n}f_{\om_i}(\varphi_{\sg^i\om}(x))\bigr).
\]         
By considering the family $F-P(F)$, we may assume without loss of generality that $P(F)=0$. Equivalently, 
$P(\tilde{f})=0$ by Proposition~3.1.4 in~\cite{gdms}, 
where 
\[ 
P(\tilde{f})=\lim_{n\rightarrow\infty}\frac{1}{n}\log Z_n(\tilde{f}) 
\]
with 
\[ 
Z_n(\tilde{f})
=\sum_{\om\in E_A^n}\exp\bigl(\sup_{\tau\in[\om]}\sum_{i=0}^{n-1}\tilde{f}(\sg^i\tau)\bigr). 
\]
Since $S$ is a CGDMS with an underlying finitely primitive matrix and $F$ is a summable H\"older family of functions, 
Theorem~3.2.3 and Proposition~4.2.5 in~\cite{gdms} assert that there exists a unique $F$-conformal measure $m_F$ 
supported on $J$. In other words, for every $\om\in E_A^*$ and for every Borel set $B\sbt X_{t(\om)}$
\[ 
m_F(\varphi_\om(B))
=\int_B\exp\bigl(S_\om(F)-P(F)|\om|\bigr)dm_F=\int_B\exp(S_\om(F))dm_F 
\]
and for all incomparable words $\om,\tau\in E_A^*$
\[ m_F\bigl(\varphi_\om(X_{t(\om)})\cap\varphi_\tau(X_{t(\tau)})\bigr)=0. \]
Recall that $S_\om(F):X_{t(\om)}\rightarrow\R$ is simply the ergodic sum
\[ 
S_\om(F)(x)=\sum_{i=1}^n f_{\om_i}(\varphi_{\sg^i\om}(x)).
\]
Moreover, $m_F=\tilde{m}_{\tilde{f}}\circ\pi^{-1}$, where $\tilde{m}_{\tilde{f}}$ is the unique eigenmeasure of the 
conjugate Perron-Frobenius operator ${\cal L}_{\tilde{f}}^*$. The existence and the uniqueness of $\tilde{m}_{\tilde{f}}$
is guaranteed by Corollary~2.7.5(a) in~\cite{gdms}. The measure $\tilde{m}_{\tilde{f}}$ is also a Gibbs state for 
$\tilde{f}$ according to Corollary~2.7.5(b). Furthermore, $\tilde{f}$ admits a unique $\sg$-invariant Gibbs state 
$\tilde{\mu}_{\tilde{f}}$ according to Corollary~2.7.5(c). This Gibbs state is completely ergodic. It is also the 
unique equilibrium measure for $\tilde{f}$ by Theorem~2.2.9. As required by that theorem, note that 
$\tilde{f}\in L^1(\tilde{\mu}_{\tilde{f}})$, for $\tilde{\zeta}\in L^1(\tilde{\mu}_{\tilde{f}})$ as shown in the 
following calculation:
\begin{equation}\label{inte}
\aligned
\int|\tilde{\zeta}|\,d\tilde{\mu}_{\tilde{f}}
   &=\sum_{e\in E}\int_{[e]}|\tilde{\zeta}|\,d\tilde{\mu}_{\tilde{f}} 
   =\sum_{e\in E}\int_{[e]}\bigl|\log|\varphi_e'(\pi(\sg\om))|\bigr|\,d\tilde{\mu}_{\tilde{f}}(\om)\\
   &\leq\sum_{e\in E}\int_{[e]}\bigl(\log K-\log\|\varphi_e'\|\bigr)d\tilde{\mu}_{\tilde{f}} \\
   &=\log K - \sum_{e\in E}\log\|\varphi_e'\|\tilde{\mu}_{\tilde{f}}([e]) \\
   &\leq\log K - C\sum_{e\in E}\log\|\varphi_e'\|\exp\bigl(\sup f_e\bigr) \\
   &\leq\log K - C\sum_{e\in E}\log\|\varphi_e'\|\exp\bigl(\sup\psi_e+u\sup\log|\varphi_e'|\bigr) \\
   &\leq\log K - Ce^{\|\Psi\|}\sum_{e\in E}\log\|\varphi_e'\|\exp\bigl(u\log K+u\log\|\varphi_e'\|\bigr) \\
   &\leq\log K - C K^u e^{\|\Psi\|}\sum_{e\in E}\|\varphi_e'\|^u\log\|\varphi_e'\| \\
   &=\log K - C K^u e^{\|\Psi\|}\sum_{e\in E}\|\varphi_e'\|^{u-\d}\,\|\varphi_e'\|^\d\log\|\varphi_e'\| \\
   &\leq\log K - Ce^{\|\Psi\|+u\log K}B\sum_{e\in E}\|\varphi_e'\|^{u-\d} \\
   &<\infty,
\endaligned
\end{equation}
where $0<\d<u-\th$ and $-\infty<B:=\inf\{\|\varphi_e'\|^\d\log\|\varphi_e'\|:e\in E\}<0$, and where we used 
the bounded distortion property (explaining the presence of $K\geq 1$), the fact that $\tilde{m}_{\tilde{f}}$ 
is a Gibbs state for $\tilde{f}$ (hence the presence of $C\geq 1$), that $P(\tilde{f})=0$, and that 
$\lim_{x\rightarrow 0^+}x^\d\log x=0$ (by means of $B$).

\noindent Since any two Gibbs states for $\tilde{f}$ are boundedly equivalent, the measures $\tilde{\mu}_{\tilde{f}}$ 
and $\tilde{m}_{\tilde{f}}$ are boundedly equivalent and thus the measures 
$\mu_F:=\tilde{\mu}_{\tilde{f}}\circ\pi^{-1}$ and $m_F=\tilde{m}_{\tilde{f}}\circ\pi^{-1}$ 
are boundedly equivalent. For this reason, $\mu_F$ is called the $S$-invariant version of $m_F$. 

Now, for every $(q,t)\in\R^2$ define the family 
\[ F_{q,t}:=qF+t\Log=q\Psi+(qu+t)\Log, \]
and its corresponding amalgamated function 
\[ \tilde{f}_{q,t}=q\tilde{f}+t\tilde{\zeta}=q\tilde{\psi}+(qu+t)\tilde{\zeta}. \]
It is easy to see that the $F_{q,t}$'s are H\"older families of functions 
of order $\g$ (the order of $F$) and that the $\tilde{f}_{q,t}$'s are H\"older 
continuous functions of order $\g$ (like $\tilde{f}$). Moreover, $F_{q,t}$ and 
$\tilde{f}_{q,t}$ are summable if and only if $qu+t>\th$.
Hence for every $(q,t)\in\R^2$ such that $qu+t>\th$, there exists a unique 
$F_{q,t}$-conformal measure $m_{F_{q,t}}$ supported on $J$. This measure is such that
$m_{F_{q,t}}=\tilde{m}_{\tilde{f}_{q,t}}\circ\pi^{-1}$, where $\tilde{m}_{\tilde{f}_{q,t}}$ is the unique 
eigenmeasure of the conjugate Perron-Frobenius operator ${\cal L}_{\tilde{f}_{q,t}}^*$. The measure 
$\tilde{m}_{\tilde{f}_{q,t}}$ is a Gibbs state for $\tilde{f}_{q,t}$. Furthermore, $\tilde{f}_{q,t}$ admits 
a unique completely ergodic $\sg$-invariant Gibbs state $\tilde{\mu}_{\tilde{f}_{q,t}}$. This Gibbs state 
is also the unique equilibrium measure for $\tilde{f}_{q,t}$. As the Gibbs states $\tilde{\mu}_{\tilde{f}_{q,t}}$ 
and $\tilde{m}_{\tilde{f}_{q,t}}$ are boundedly equivalent, so are the measures 
$\mu_{F_{q,t}}:=\tilde{\mu}_{\tilde{f}_{q,t}}\circ\pi^{-1}$ and 
$m_{F_{q,t}}=\tilde{m}_{\tilde{f}_{q,t}}\circ\pi^{-1}$. 
Note that $\tilde{f}_{q,t}\in L^1(\tilde{\mu}_{q,t})$, for $\tilde{\zeta}\in L^1(\tilde{\mu}_{q,t})$ 
as a calculation similar to~(\ref{inte}) shows.
Furthermore, 
\begin{equation}\label{press}
P\bigl((qu+t)\Log\bigr)-|q|\,\|\Psi\| 
   \leq P(F_{q,t})\leq 
   P\bigl((qu+t)\Log\bigr)+|q|\,\|\Psi\|,
\end{equation}
where $P((qu+t)\Log)=P(qu+t)$ (cf. section~\ref{prelgdms}).
We now state fundamental properties of the pressure as a function of the two variables
$q$ and $t$.

\bthm\lab{1mam1}
Let $S=\{\varphi_e:X_{t(e)}\rightarrow X_{i(e)}\,|\,e\in E\}$ be a cofinitely 
regular CGDMS satisfying the Strong Open Set Condition (SOSC) and having an 
underlying finitely primitive edge incidence matrix $A$. 
Let $\th$ be the finiteness parameter of $S$, and let $u>\th$. Then
the pressure function $(q,t)\mapsto 
P(q,t):=P(F_{q,t})=P(\tilde{f}_{q,t})$, $(q,t)\in\R^2$ satisfies the following properties.
\begin{itemize}
\item[(a)] $P(q,t)<\infty$ if and only if $qu+t>\th$; 
\item[(b)] 
If 
$ 
(q_2-q_1)(\sup\Psi+u\log s)+(t_2-t_1)\log s\leq 0, 
$
then
$
P(q_2,t_2)\leq P(q_1,t_1). 
$
\noindent In particular, if $\sup\Psi\leq-u\log s$, then $P(q,t)$ is decreasing with respect to both variables 
$q\in\R$ and $t\in\R$. 
\item[(c)] $t\mapsto P(q,t)$ is strictly decreasing on $(\th-qu,\infty)$;
\item[(d)] $\lim_{t\to\infty}P(q,t)=-\infty$;
\item[(e)] $\lim_{t\to(\th-qu)^+}P(q,t)=\infty$;
\item[(f)] $\frac{\partial P}{\partial t}(q,t)=-\chi_{\tilde{\mu}_{q,t}}$ for every $(q,t)\in\R^2$ 
           such that $qu+t>\th$, where 
$\chi_{\tilde{\mu}_{q,t}}:=-\int \tilde{\zeta}\,d\tilde{\mu}_{q,t}$ is the Lyapunov exponent of 
$\tilde{\mu}_{q,t}$;
\item[(g)] $t\mapsto P(q,t)$ is convex (and thereby continuous) on the interval $(\th-qu,\infty)$;
\end{itemize}
\ethm  

{\sl Proof.} (a) This follows from~(\ref{press}) and the fact that 
$P(\overline{t}):=P(\overline{t}\,\Log)<\infty$ if and only if $\overline{t}>\th$.

(b) Let $q_1\leq q_2$, $t_1\leq t_2$. If $q_1u+t_1\leq\th$, then $P(q_1,t_1)=\infty$ and the statement thus holds. 
So suppose that $q_1u+t_1>\th$. Then the $n$-th partition function of $\tilde{f}_{q_2,t_2}$ satisfies
\begin{eqnarray*}
Z_n(\tilde{f}_{q_2,t_2})
              &=&\sum_{\om\in E_A^n}\exp\bigl(\sup_{\tau\in[\om]}S_n\tilde{f}_{q_2,t_2}(\tau)\bigr) \\
              &=&\sum_{\om\in E_A^n}\sup_{\tau\in[\om]}\exp\bigl(S_n\tilde{f}_{q_2,t_2}(\tau)\bigr) \\
              &=&\sum_{\om\in E_A^n}\sup_{\rho\in E_A^\infty:A_{\om_n\rho_1}=1}
              \Bigl(\exp\bigl(S_n(q_2\tilde{\psi}(\om\rho))\bigr)|\varphi_\om'(\pi(\rho))|^{q_2u+t_2}\Bigr) \\
              &\le&\sum_{\om\in E_A^n}\sup_{\tau\in[\om]}\Bigl(\exp\bigl(S_n(q_2\tilde{\psi})(\tau)\bigr)\Bigr)
              \|\varphi_\om'\|^{q_2u+t_2} \\
              &\le&\sum_{\om\in E_A^n}\sup_{\tau\in[\om]}\Bigl(\exp\bigl(S_n(q_1\tilde{\psi})(\tau)
                                                            + S_n((q_2-q_1)\tilde{\psi})(\tau)\bigr)\Bigr)
              \|\varphi_\om'\|^{q_1u+t_1}\|\varphi_\om'\|^{(q_2-q_1)u+(t_2-t_1)} \\
              &\le&s^{n[(q_2-q_1)u+(t_2-t_1)]}K^{q_1u+t_1}e^{n(q_2-q_1)\sup\Psi}
                   \sum_{\om\in E_A^n}\sup_{\tau\in[\om]}\exp\bigl(S_n(q_1\tilde{\psi})(\tau)\bigr)
                                                   \inf_{x\in X_{t(\om)}}|\varphi_\om'(x)|^{q_1u+t_1} \\
              &\le&e^{n(q_2-q_1)\sup\Psi}s^{n[(q_2-q_1)u+(t_2-t_1)]}K^{q_1u+t_1}
                   \sum_{\om\in E_A^n}\sup_{\rho\in E_A^\infty:A_{\om_n\rho_1}=1}
                          \Bigl(\exp\bigl(S_n(q_1\tilde{\psi})(\om\rho)\bigr)
                          |\varphi_\om'(\pi(\rho))|^{q_1u+t_1}\Bigr) \\
              &=&e^{n(q_2-q_1)\sup\Psi}s^{n[(q_2-q_1)u+(t_2-t_1)]}K^{q_1u+t_1}Z_n(\tilde{f}_{q_1,t_1}).
\end{eqnarray*}
Therefore 
\begin{equation}\label{estp}
\aligned
P(q_2,t_2)
&=\lim_{n\rightarrow\infty}\frac{1}{n}\log Z_n(\tilde{f}_{q,t_2}) \\
&\le(q_2-q_1)\sup\Psi+[(q_2-q_1)u+(t_2-t_1)]\log s
     +\lim_{n\rightarrow\infty}\frac{1}{n}\log Z_n(\tilde{f}_{q_1,t_1}) \\
&=(q_2-q_1)\sup\Psi+[(q_2-q_1)u+(t_2-t_1)]\log s+P(q_1,t_1).
\endaligned
\end{equation}
Part (b) follows immediately.

\

(c) Letting $q_1=q_2=q$ and $t_1<t_2$ in~(\ref{estp}) gives (c).

\

(d) This also follows from~(\ref{estp}) by setting $q_1=q_2=q$, 
$t_1>\th$ and $t_2=t$ and letting $t\rightarrow\infty$.
  
\

(e) Let $t>\th-qu$. 
Then
\begin{eqnarray*}
Z_n(\tilde{f}_{q,t})&\ge&\sum_{\om\in E_A^n}\exp(-nq\inf\Psi)K^{-(qu+t)}\|\varphi_\om'\|^{qu+t} \\
              &=&\exp(-nq\inf\Psi)K^{-(qu+t)}\sum_{\om\in E_A^n}\|\varphi_\om'\|^{qu+t}.
\end{eqnarray*}
Therefore 
\[ 
P(q,t)\ge-q\inf\Psi+P\bigl((qu+t)\Log\bigr)=-q\inf\Psi+P(qu+t). 
\]
Thus, 
\[ 
\lim_{t\rightarrow(\th-qu)^+}P(q,t)
\ge-q\inf\Psi+\lim_{t\rightarrow(\th-qu)^+}P(qu+t)
=-q\inf\Psi+P(\th)
=\infty 
\]
since $S$ is cofinitely regular.

\

(f) This follows from Proposition~2.6.13 in~\cite{gdms} and the fact that 
$-\chi_{\tilde{\mu}_{q,t}}=\int \tilde{\zeta}\,d\tilde{\mu}_{q,t}>-\infty$ 
whenever $qu+t>\th$ by a calculation similar to~(\ref{inte}).
 
\

(g) This follows immediately from Proposition~2.6.14 in~\cite{gdms}.
\endpf

\

Observe that in the proof of part (e) of the above theorem $P(\th)\ge q\inf\Psi$ suffices to guarantee the
existence of a zero for the pressure function $t\mapsto P(q,t)$. The assumption of cofinite regularity 
on $S$ ensures that the pressure function $t\mapsto P(q,t)$ has a zero for every $q\in\R$. 

\bcor\lab{c1mam5}
For all $q\in\R$ there exists a unique $T(q)\in(\th-qu,\infty)$ such that 
$P(q,T(q))=0$. The function $T(q)$ is called the 
temperature function.
\ecor

\

In order to allege notation, let 
\[ 
\tilde{f}_q=\tilde{f}_{q,T(q)},\,\,\,\,\, F_q=F_{q,T(q)},\,\,\,\,\, 
\tilde{m}_q=\tilde{m}_{\tilde{f}_{q,T(q)}},\,\,\,\,\, m_q=m_{F_{q,T(q)}},\,\,\,\,\, 
\tilde{\mu}_q=\tilde{\mu}_{\tilde{f}_{q,T(q)}},\,\,\,\,\, \mu_q=\mu_{F_{q,T(q)}}. 
\]

\


Now, let $q\in\R$. Recall that $\tilde{\zeta}\in L^1(\tilde{\mu}_q)$ and thus $\tilde{f}\in L^1(\tilde{\mu}_q)$.
Thereafter, let 
\[ 
\a(q)=\frac{\int\tilde{f}\,d\tilde{\mu}_q}{-\chi_{\tilde{\mu}_q}(\sg)}
     =\frac{\int\tilde{f}\,d\tilde{\mu}_q}{\int\tilde{\zeta}\,d\tilde{\mu}_q}.
\]
By the variational principle for pressure (cf. Theorems~2.1.6--2.1.8 in~\cite{gdms}), note that  
\begin{equation}\label{intf}
\int\tilde{f}\,d\tilde{\mu}_q
           \leq P(\tilde{f})-h_{\tilde{\mu}_q}(\sg)
           =-h_{\tilde{\mu}_q}(\sg)
           \leq 0. 
\end{equation}
Hence $\a(q)\geq 0$. Moreover,
\[ \a(q)\leq\frac{\int|\tilde{\psi}+u\tilde{\zeta}|\,d\tilde{\mu}_q}{\int|\tilde{\zeta}|\,d\tilde{\mu}_q}
        \leq u+\frac{\int|\tilde{\psi}|\,d\tilde{\mu}_q}{\int|\tilde{\zeta}|\,d\tilde{\mu}_q}
        \leq u+\frac{\|\Psi\|}{-\log s}. \]
Thus, $0\leq\a(q)<\infty$.

\

Finally, we study some basic properties of the temperature function $T(q)$.

\bthm\label{alphat}
The temperature function $q\mapsto T(q)$ exhibits the following properties.
\begin{itemize}
\item[(a)] The function $T:\R\rightarrow\R$ is real-analytic;
\item[(b)] $T(0)=\HD(J)$ while $T(1)=0$;
\item[(c)] $T'(q)=-\a(q)<0$ for all $q\in\R$;
\item[(d)] The function $q\mapsto T(q)$, $q\in \R$, is convex, meaning that 
$T''(q)\geq 0$ for all $q\in\R$. This function is not
strictly convex if and only if $\tilde{\mu}_{\tilde f}$ is equal to 
$\tilde{\mu}_{-\HD(J)\tilde\zeta}$.
\end{itemize} 
\ethm

{\sl Proof.} (a) By Proposition~2.6.13 in~\cite{gdms}, 
$\frac{\partial P}{\partial t}(q,t)=\int\tilde{\zeta}\,d\tilde{\mu}_{q,t}
=-\chi_{\tilde{\mu}_{q,t}}(\sg)<0$ 
for every $(q,t)\in\R^2$ such that $qu+t>\th$. In particular, this is true for all pairs $(q,T(q))$. 
Since $T(q)$ is uniquely determined by the condition
$P(q,T(q))=0$, it follows from Theorem~2.6.12 in~\cite{gdms} and the implicit function theorem that
$T$ is real-analytic on $\R$.

\

(b) Since $S$ is regular, we have $P(\HD(J)\Log)=0$, which means that $T(0)=\HD(J)$. 
Moreover, since $\tilde{f}_{1,0}=\tilde{f}$ and $P(F)=P(\tilde{f})=0$ by assumption, 
we deduce that $P(1,0)=P(\tilde{f}_{1,0})=P(\tilde{f})=0$, and from the uniqueness of $T(1)$ 
it follows that $T(1)=0$.

\

(c) It follows from the fact that $P(q,T(q))=0$ for all $q\in\R$ 
and from Proposition~2.6.13 in~\cite{gdms} that 
\[ 0=\frac{dP}{dq}(q,T(q))=\frac{\partial P}{\partial q}(q,T(q))
                           +\frac{\partial P}{\partial t}(q,T(q))\cdot T'(q)
                          =\int \tilde{f}\,d\tilde{\mu}_q
                          -\chi_{\tilde{\mu}_q}(\sg)T'(q).       
\]
Hence
\[ 
T'(q)
=-\frac{\int\tilde{f}\,d\tilde{\mu}_q}{-\chi_{\tilde{\mu}_q}(\sg)}
=-\a(q). 
\]
Having already observed that $\a(q)\geq 0$, we thus know that $T'(q)\leq 0$. In order to 
prove that $T'(q)<0$, we need to show that $\int\tilde{f}\,d\tilde{\mu}_q\neq 0$. But since 
$\int\tilde{f}\,d\tilde{\mu}_q\leq-h_{\tilde{\mu}_q}(\sg)\leq 0$ by~(\ref{intf}), it suffices to show that
$h_{\tilde{\mu}_q}(\sg)>0$. This follows immediately from Theorem~2.5.2 in~\cite{gdms}. 
\

(d) Lemma~4.9.5 in~\cite{gdms} (with $\Delta_1$ replaced by $\R$) 
and its proof carry over to the current setting without any change.  
\endpf

\section{Multifractal analysis of the conformal measure $m_F$ over a subset of $J$}\label{jr}

Let $S=\{\varphi_e:X_{t(e)}\rightarrow X_{i(e)}\,|\,e\in E\}$ be a cofinitely regular CGDMS 
satisfying SOSC 
and having an underlying edge incidence matrix 
$A$ which is finitely primitive. Let $\th$ be the finiteness parameter of $S$. 
Let also $F$ be a family of functions of the form $F=\Psi+u\Log$ such that $P(F)=0$,
where $\Psi$ is a bounded H\"older family of functions and $u>\th$.  

\

We shall now develop the multifractal analysis of the conformal measure $m_F$ 
(or equivalently, the invariant measure $\mu_F$) associated to the family $F$. We shall 
conduct this analysis by means of balls and we shall restrict ourselves to a subset $J_r$ 
of the limit set $J$ of $S$.  
As we shall see later, $J_r$ is often a fairly large subset of $J$. 
By definition, $J_r$ is the set of points of $J$ which are coded by 
the set of infinite admissible words 
\[ 
E_r^\infty=\Bigl\{\om\in E_A^\infty\,\,\bigl|\,\,\limsup_{n\rightarrow\infty}
                                                \dist\bigl(\pi(\sg^n\om),\partial X_{i(\sg^n\om)}\bigr)>0\Bigr\}. 
\] 
The words in this set code points of the limit set that behave tamely when a multifractal analysis 
is carried out using balls, for infinitely many of their iterates are positively separated from the 
boundary of the phase space. The conformality of the measure $m_F$ can then be used at those 
iterates to estimate the local dimension of $m_F$ from above. Before going any further, 
we observe that $E_r^\infty$ is a set of full measure. 


\blem
For every ergodic, $\sg$-invariant Borel probability measure $\tilde{\mu}$ on $E_A^\infty$ with 
$\supp\tilde{\mu}=E_A^\infty$, we have $\tilde{\mu}(E_r^\infty)=1$.
\elem

{\sl Proof.} Let $\tilde{\mu}$ be an ergodic, $\sg$-invariant Borel probability measure on 
$E_A^\infty$ with $\supp\tilde{\mu}=E_A^\infty$. Observe that $E_r^\infty$ is completely 
$\sg$-invariant, that is, $\sg^{-1}(E_r^\infty)=E_r^\infty=\sg(E_r^\infty)$. Thus, 
by ergodicity of $\tilde{\mu}$, we have that $\tilde{\mu}(E_r^\infty)$ is $0$ or $1$. We 
shall now show that this latter possibility always prevails. Since $S$ satisfies SOSC, 
there exists $x\in J_v\cap\Int(X_v)$ for some $v\in V$. Let $\om\in E_A^\infty$ be such 
that $\pi(\om)=x$. Let also $0<r<\dist(x,\partial X_v)$. Note that for any $\tau\in E_A^*$, we have 
$\pi([\tau])\sbt\varphi_\tau(X_{t(\tau)})$. Since 
$\varphi_{\om|_k}(X_{t(\om_k)})\sbt B(\pi(\om),r)=B(x,r)$ for all $k\in\N$ large enough,
we obtain that $[\om|_k]\sbt\pi^{-1}(B(x,r))$ for some $k\in\N$. 
Then $\tilde{\mu}(\pi^{-1}(B(x,r)))\ge\tilde{\mu}([\om|_k])>0$
since $\supp\tilde{\mu}=E_A^\infty$. It follows from Birkhoff's Ergodic Theorem 
that the set of infinite admissible words whose iterates' images visit infinitely 
many times the ball $B(x,r)$ has measure $1$, that is,
\[ \tilde{\mu}\Bigl(\bigl\{\tau\in E_A^\infty\,\bigl|\,\sg^n\tau\in\pi^{-1}(B(x,r)) 
                           \mbox{ for infinitely many } n\mbox{'s}\bigr\}\Bigr)=1. \] 
Therefore $\tilde{\mu}(E_r^\infty)=1$. (The same conclusion can be drawn by means of 
Poincar\'e's Recurrence Theorem.) 
\endpf

\

This lemma tells us that 
$\tilde{\mu}_{\tilde{f}}(E_r^\infty)=1$ and 
$\tilde{\mu}_{\tilde{f}_{q,t}}(E_r^\infty)=1$ 
for all $(q,t)\in\R^2$ such that $qu+t>\th$.
In particular, $\tilde{\mu}_q(E_r^\infty)=1$ 
for all $q\in\R$.
This implies immediately that 
$\mu_F(J_r)=\tilde{\mu}_{\tilde{f}}\circ\pi^{-1}(\pi(E_r^\infty))=1$,
that is, the set $E_r^\infty$ is a set of full $\mu_F$-measure. Similarly,
$\mu_{F_{q,t}}(J_r)=1$ for all $(q,t)\in\R^2$ such that $qu+t>\th$.
In particular, $\mu_q(J_r)=1$ for all $q\in\R$.

\

In fact, $J_r$ contains a rich family of subsets of full measure. 
To define these subsets, we proceed as follows. 
For every $\om\in E_A^\infty$ and $r\geq 0$, let $\{n_j(\om,r)\}$ be the increasing sequence 
of all positive integers $n$ such that $\sg^n\om\in\pi^{-1}(X_{i(\sg^n\om)}\backslash 
B(\partial X_{i(\sg^n\om)},r))$. This sequence may be empty, non-empty and finite, or 
infinite depending on $\om$ and $r$. However, for every $\om\in E_r^\infty$ the sequence 
$\{n_j(\om,r)\}$ is infinite for every $0\leq r<r_{max}(\om)$, where
\[ r_{max}(\om):=\limsup_{n\rightarrow\infty}\dist\bigl(\pi(\sg^n\om),\partial X_{i(\sg^n\om)}\bigr)>0. \]  

\

\noindent Now, for every $R\geq 0$ define the completely invariant set 
\[ 
E_{rr}^\infty(R):= \Bigl\{\om\in E_r^\infty\,|\,R<r_{max}(\om)
\mbox{ and }  
\lim_{j\rightarrow\infty} \frac{S_{n_{j+1}(\om,R)}\tilde{\zeta}(\om)-S_{n_j(\om,R)}\tilde{\zeta}(\om)}
                               {S_{n_j(\om,R)}\tilde{\zeta}(\om)}=0
\Bigr\}. 
\]
We claim that these subsets of $E_r^\infty$ have all full measure.
\blem\lab{l1mam11}
For every $R\geq 0$ small enough, we have that $\tilde{\mu}(E_{rr}^\infty(R))=1$ for all ergodic, 
$\sg$-invariant Borel probability measure $\tilde{\mu}$ on $E_A^\infty$ with 
$\supp\tilde{\mu}=E_A^\infty$ and such that $\int \tilde{\zeta} \, d\tilde{\mu} > -\infty$. 
\elem 

{\sl Proof.} Let $0\leq R<D$, where $D=\sup_{x\in J}\dist(x,\partial X)
:=\sup_{v\in V}\sup_{x\in J_v}\dist(x,\partial X_v)$. SOSC guarantees that $D>0$. 
Let $\tilde{\mu}$ be an ergodic, $\sg$-invariant Borel probability measure on $E_A^\infty$ 
with $\supp\tilde{\mu}=E_A^\infty$ and such that $\int \tilde{\zeta} \, d\tilde{\mu} > -\infty$. 
Since $E_{rr}^\infty(R)$ is completely $\sg$-invariant, the ergodicity of $\tilde{\mu}$ forces 
$\tilde{\mu}(E_{rr}^\infty(R))$ to equal $0$ or $1$. We shall now prove that this latter possibility 
always prevails. Since $0\leq R<D$, there is $x\in J_v\backslash B(\partial X_v,(R+D)/2)$ for 
some $v\in V$. Let 
$\om\in E_A^\infty$ be such that $\pi(\om)=x$. Then there is some $k\in\N$ 
such that $\pi([\om|_k])\sbt\varphi_{\om|_k}(X_{t(\om_k)})\sbt B(\pi(\om),(D-R)/2)=B(x,(D-R)/2)$, 
or equivalently $[\om|_k]\sbt\pi^{-1}B(x,(D-R)/2)$. Then 
\[ 
\tilde{\mu}\bigl(\pi^{-1}(X_v\backslash B(\partial X_v,R))\bigr)
\geq\tilde{\mu}\bigl(\pi^{-1}B(x,(D-R)/2)\bigr)\ge\tilde{\mu}([\om|_k])>0 
\] 
since $\supp\tilde{\mu}=E_A^\infty$. Applying Birkhoff's Ergodic Theorem twice (once with the 
characteristic function of the set $\pi^{-1}(X_v\backslash B(\partial X_v,R))$ and once with the 
potential function $\tilde{\zeta}$), we obtain that 
the set 
\[ 
\Bigl\{\om\in E_r^\infty\,|\,R<r_{max}(\om),\, \lim_{j\rightarrow\infty}\frac{n_{j+1}(\om,R)}{n_j(\om,R)}=1 
                           \,\mbox{ and }\, 
                           \lim_{j\rightarrow\infty}\frac{S_{n_j(\om,R)}\tilde{\zeta}(\om)}{n_j(\om,R)}
                                      =\int \tilde{\zeta}\, d\tilde{\mu}=-\chi_{\tilde{\mu}}(\sg) \Bigr\} 
\]
has measure $1$. Writing $n_j$ instead of $n_j(\om,R)$ to allege notation, 
we have for every $\om$ in this set that 
\begin{eqnarray*}
\lim_{j\rightarrow\infty} \frac{S_{n_{j+1}}\tilde{\zeta}(\om)-S_{n_j}\tilde{\zeta}(\om)}
                               {S_{n_j}\tilde{\zeta}(\om)}
&=&\lim_{j\rightarrow\infty} \frac{n_{j+1}}{n_j}\frac{\frac{1}{n_{j+1}}\bigl(S_{n_{j+1}}\tilde{\zeta}(\om)
                                                                            -S_{n_j}\tilde{\zeta}(\om)\bigr)
                                                     }
                                                     {\frac{1}{n_j}S_{n_j}\tilde{\zeta}(\om)} \\
&=&\lim_{j\rightarrow\infty} \frac{\frac{1}{n_{j+1}}S_{n_{j+1}}\tilde{\zeta}(\om)
                                   -\frac{n_j}{n_{j+1}}\cdot\frac{1}{n_{j}}S_{n_j}\tilde{\zeta}(\om)}
                                  {\frac{1}{n_j}S_{n_j}\tilde{\zeta}(\om)} \\
&=&\frac{-\chi_{\tilde{\mu}}(\sg)-1\cdot(-\chi_{\tilde{\mu}}(\sg))}{-\chi_{\tilde{\mu}}(\sg)}=0. 
\end{eqnarray*}
It follows immediately that $\tilde{\mu}(E_{rr}^\infty(R))=1$.
\endpf

\

This lemma reveals that 
$\tilde{\mu}_{\tilde{f}}(E_{rr}^\infty(R))=1$ and 
$\tilde{\mu}_{\tilde{f}_{q,t}}(E_{rr}^\infty(R))=1$ 
for all $(q,t)\in\R^2$ such that $qu+t>\th$.
In particular, $\tilde{\mu}_q(E_{rr}^\infty(R))=1$ 
for all $q\in\R$.
This implies immediately that 
$\mu_F(\pi(E_{rr}^\infty(R)))=1$,
that is, the set $E_{rr}^\infty(R)$ is a set of full $\mu_F$-measure for all $R\in\R$. Consequently,
$\mu_{F_{q,t}}(\pi(E_{rr}^\infty(R)))=1$ for all $(q,t)\in\R^2$ such that $qu+t>\th$.
In particular, $\mu_q(\pi(E_{rr}^\infty(R)))=1$ for all $q\in\R$.

\

As an immediate corollary, we obtain that the completely invariant set 
\[ E_{rr}^\infty=\bigcup_{R>0}E_{rr}^\infty(R) \]
is a set of full measure.


\bcor\lab{l1mam111}
For every ergodic, $\sg$-invariant Borel probability measure $\tilde{\mu}$ on $E_A^\infty$ 
with $\supp\tilde{\mu}=E_A^\infty$ and such that $\int \tilde{\zeta} \, d\tilde{\mu} > -\infty$, 
we have $\tilde{\mu}(E_{rr}^\infty)=1$.
\ecor 

%
%
%
%

\

We now recall a few basic definitions from multifractal analysis.

\

Let $\mu$ be a Borel probability measure on $X$. The pointwise or local 
dimension $d_\mu(x)$ of $\mu$ at $x\in X$  
is the power law behaviour (if any) of $\mu(B(x,r))$ for small $r>0$, that is, 
\[ d_\mu(x)=\lim_{r\rightarrow 0}\frac{\log\mu(B(x,r))}{\log r}. \]
We further define the lower and upper dimensions of $\mu$ at $x\in X$ by
\[ \underline{d}_\mu(x)=\liminf_{n\rightarrow\infty}\frac{\log\mu(B(x,r))}{\log r}\]
and
\[ \overline{d}_\mu(x)=\limsup_{n\rightarrow\infty}\frac{\log\mu(B(x,r))}{\log r},\]
%
respectively. Denote the set of points of $J_r$ at which the local dimension of a measure $\mu$ is equal to $\a$ by
\[ J_{r,\mu}(\a)=\{x\in J_r\,|\,d_\mu(x)=\a\}. \]
Denote the Hausforff dimension of $J_{r,\mu}(\a)$ by 
\[ f_{r,\mu}(\a)=\HD(J_{r,\mu}(\a)). \]  
Now for every $\a\geq 0$, let 
\[ 
E_r^\infty(\a)=\bigl\{\om\in E_r^\infty:
\lim_{n\rightarrow\infty}\frac{S_n\tilde{f}(\om)}{S_n\tilde{\zeta}(\om)}=\a\bigr\}. 
\]

\

We shall now prove that for every $q\in\R$ the measure $\mu_q$ confers full 
measure to the set of points of $J_r$ where the local dimension of the measure 
$m_F$ is $\a(q)$.

\bthm\lab{p2mam11}
The following statements hold.
\begin{itemize}
\item[(a)] For every $\a\ge 0$, we have $\pi(E_r^\infty(\a)\cap E_{rr}^\infty)\sbt J_{r,m_F}(\a)$;
\item[(b)] $\tilde{\mu}_q(E_r^\infty(\a(q))\cap E_{rr}^\infty)=1$ for all $q\in\R$;
\item[(c)] $\mu_q(J_{r,m_F}(\a(q)))=1$ for all $q\in\R$.
\end{itemize}
\ethm

{\sl Proof.} (a) Let $x\in \pi(E_r^\infty(\a)\cap E_{rr}^\infty)$. Then there is some 
$\om\in E_r^\infty(\a)\cap E_{rr}^\infty$ such that $\pi(\om)=x$. Therefore 
$\om\in E_{rr}^\infty(R)$ for some $0<R<\min_{v\in V}\dist(X_v,\partial W_v)$. Let 
$\{n_j\}_{j\in\N}:=\{n_j(\om,R)\}_{j\in\N}$ be the increasing sequence of all $n$'s 
such that $\pi(\sg^n\om)\in X_{t(\om_n)}\backslash B(\partial X_{t(\om_n)},R)$. Let 
$0<r\leq K^{-1}R|\varphi_{\om|_{n_1}}'(\pi(\sg^{n_1}\om))|$. Let $j\in\N$ be the 
unique natural number so that 
\[ 
K^{-1}R|\varphi_{\om|_{n_{j+1}}}'(\pi(\sg^{n_{j+1}}\om))|
<r\leq K^{-1}R|\varphi_{\om|_{n_j}}'(\pi(\sg^{n_j}\om))|. 
\]

\noindent Since $B(\pi(\sg^{n_j}\om),R)\sbt W_{t(\om_{n_j})}$, the conformality of 
the generators of the system ensures that
\[ 
B(x,r)\sbt B\bigl(\pi(\om),K^{-1}R|\varphi_{\om|_{n_j}}'(\pi(\sg^{n_j}\om))|\bigr)
      \sbt \varphi_{\om|_{n_j}}\bigl(B(\pi(\sg^{n_j}\om),R)\bigr), 
\]
where the last inclusion follows from relation (4.22) in~\cite{gdms}. 
%
%
%
%
%
%
%
%
%
Moreover, every $y\in J_{t(\om_{n_j})}$ admits a $\tau\in [\om|_{n_j}]$ such that $y=\pi(\sg^{n_j}\tau)$, 
and for such $y$ and $\tau$ we have $S_{\om|_{n_j}}F(y)=S_{n_j}\tilde{f}(\tau)$. 
The conformality of $m_F$, the fact that $B(\pi(\sg^{n_j}\om),R)\sbt\Int(X_{t(\om_{n_j})})$, 
and the OSC then give
\begin{equation}\lab{1mam14}
\aligned
m_F(B(x,r))
&\leq m_F\bigl(\varphi_{\om|_{n_j}}(B(\pi(\sg^{n_j}\om),R))\bigr) \\
&\leq m_F\bigl(\varphi_{\om|_{n_j}}(\Int(X_{t(\om_{n_j})}))\bigr) \\
&=    m_F\bigl(\varphi_{\om|_{n_j}}(\Int(X_{t(\om_{n_j})}))\cap J_{i(\om)}\bigr) \\
&=    m_F\Bigl(\varphi_{\om|_{n_j}}\bigl(\Int(X_{t(\om_{n_j})})\cap J_{t(\om_{n_j})}\bigr)\Bigr) \\
&\leq \exp\bigl(\sup_{y\in \Int(X_{t(\om_{n_j})})\cap J_{t(\om_{n_j})}} S_{\om|_{n_j}}F(y)\bigr) 
      \ m_F\bigl(\Int(X_{t(\om_{n_j})})\cap J_{t(\om_{n_j})}\bigr) \\
&\leq\exp\bigl(\sup_{\tau\in[\om|_{n_j}]} S_{n_j}\tilde{f}(\tau)\bigr) \\
&\leq B(\tilde{f}) \exp\bigl(S_{n_j}\tilde{f}(\om)\bigr),
\endaligned
\end{equation}
where $B(\tilde{f})$ is a constant of bounded variation for $\tilde{f}$ (see the bounded variation 
principle for H\"older continuous potentials, Lemma~2.3.1 in~\cite{gdms}). 
%
%

\

\noindent On the other hand, the conformality of the generators of the system guarantees that 

\[ 
B(x,r)\spt B\bigl(\pi(\om),K^{-1}R|\varphi_{\om|_{n_{j+1}}}'(\pi(\sg^{n_{j+1}}\om))|\bigr) 
      \spt \varphi_{\om|_{n_{j+1}}}\bigl(B(\pi(\sg^{n_{j+1}}\om),K^{-2}R)\bigr).
\]

\noindent Moreover, every $y\in J_{t(\om_{n_{j+1}})}$ admits a $\tau\in [\om|_{n_{j+1}}]$ 
such that $y=\pi(\sg^{n_{j+1}}\tau)$, and for such $y$ and $\tau$ we have 
$S_{\om|_{n_{j+1}}}F(y)=S_{n_{j+1}}\tilde{f}(\tau)$. 
Then the conformality of the measure $m_F$ leads to 
\begin{equation}\label{1mam15}
\aligned  
m_F(B(x,r))
&\geq m_F\bigl(\varphi_{\om|_{n_{j+1}}}(B(\pi(\sg^{n_{j+1}}\om),K^{-2}R))\bigr) \\
&\geq m_F\Bigl(\varphi_{\om|_{n_{j+1}}}\bigl(B(\pi(\sg^{n_{j+1}}\om),K^{-2}R)
               \cap J_{t(\om_{n_{j+1}})}\bigr)
         \Bigr) \\
&\geq \exp\bigl(\inf_{y\in  J_{t(\om_{n_{j+1}})}} S_{\om|_{n_{j+1}}}F(y)\bigr) 
                \ m_F\bigl(B(\pi(\sg^{n_{j+1}}\om),K^{-2}R)\cap J_{t(\om_{n_{j+1}})}\bigr) \\
&\geq \exp\bigl(\inf_{\tau\in[\om|_{n_{j+1}}]} S_{n_{j+1}}\tilde{f}(\tau)\bigr) 
                \ m_F\bigl(B(\pi(\sg^{n_{j+1}}\om),K^{-2}R)\bigr) \\
&\geq B(\tilde{f})\exp(S_{n_{j+1}}\tilde{f}(\om))\ M_F(K^{-2}R) \\
&=M_F(K^{-2}R) B(\tilde{f})\ \exp(S_{n_{j+1}}\tilde{f}(\om)),  
\endaligned
\end{equation}
where 
\[ M_F(a)=\inf_{y\in\overline{J}}m_F(B(y,a))>0 \]
and $B(\tilde{f})$ is a constant of bounded variation for $\tilde{f}$ (see the bounded 
variation principle for H\"older continuous potentials, Lemma~2.3.1 in~\cite{gdms}). 

\

\noindent From the definition of $n_j$, we also have that
\begin{equation}\label{1mam16}
\log(K^{-1}R)+\log|\varphi_{\om|_{n_{j+1}}}'(\pi(\sg^{n_{j+1}}\om))|
<\log r\le\log(K^{-1}R)+\log|\varphi_{\om|_{n_j}}'(\pi(\sg^{n_j}\om))|<0. 
\end{equation}

\noindent As $\mu_F$ and $m_F$ are boundedly equivalent, we deduce from~(\ref{1mam15}) 
and~(\ref{1mam16}) that
\begin{eqnarray*}
d_{m_F}(x)
=\lim_{r\rightarrow 0}\frac{\log m_F(B(x,r))}{\log r} 
&\le&\lim_{j\rightarrow\infty}\frac{S_{n_{j+1}}\tilde{f}(\om)+\log(M_F(K^{-2}R) B(\tilde{f}))}
                                   {\log(K^{-1}R)+\log|\varphi_{\om|_{n_j}}'(\pi(\sg^{n_j}\om))|} \\
&=&\lim_{j\rightarrow\infty}\frac{S_{n_{j+1}}\tilde{f}(\om)+\log(M_F(K^{-2}R) B(\tilde{f}))}
                                 {\log(K^{-1}R)+S_{n_j}\tilde{\zeta}(\om)} \\
&=&\lim_{j\rightarrow\infty}\frac{\frac{S_{n_{j+1}}\tilde{f}(\om)}{S_{n_{j+1}}\tilde{\zeta}(\om)}
                                 +\frac{\log(M_F(K^{-2}R) B(\tilde{f}))}{S_{n_{j+1}}\tilde{\zeta}(\om)}
                                 }
                                 {\frac{\log(K^{-1}R)}{S_{n_{j+1}}\tilde{\zeta}(\om)}
                                 +\frac{S_{n_j}\tilde{\zeta}(\om)}{S_{n_{j+1}}\tilde{\zeta}(\om)}
                                 } \\
&=&\lim_{j\rightarrow\infty}\frac{\frac{S_{n_{j+1}}\tilde{f}(\om)}{S_{n_{j+1}}\tilde{\zeta}(\om)}
                                 +\frac{\log(M_F(K^{-2}R)B(\tilde{f}))}{S_{n_{j+1}}\tilde{\zeta}(\om)}
                                 }
                                 {\frac{\log(K^{-1}R)}{S_{n_{j+1}}\tilde{\zeta}(\om)}
                                 +\Bigl(1+\frac{(S_{n_{j+1}}-S_{n_j})\tilde{\zeta}(\om)}
                                               {S_{n_j}\tilde{\zeta}(\om)}\Bigr)^{-1}
                                 } \\
&=&\frac{\a+0}{0+(1+0)^{-1}}=\a
\end{eqnarray*}
since $\om\in E_r^\infty(\a)\cap E_{rr}^\infty(R)$.

\

\noindent Similarly, we deduce from~(\ref{1mam14}) and~(\ref{1mam16}) that
\[
d_{m_F}(x)
\ge\lim_{j\rightarrow\infty}\frac{S_{n_j}\tilde{f}(\om)+\log B(\tilde{f})}
                                 {\log(K^{-1}R)+\log|\varphi_{\om|_{n_{j+1}}}'(\pi(\sg^{n_{j+1}}\om))|} 
=\a. 
\]

\


\noindent Hence $d_{m_F}(x)=\a$. 
This completes the proof of (a).

\

(b) Let $q\in\R$. According to Birkhoff's Ergodic Theorem and Corollary~\ref{l1mam111}, there exists $E_q\sbt E_{rr}^\infty$ 
such that $\tilde{\mu}_q(E_q)=1$ and so that for all $\om\in E_q$ we have
\[ 
\lim_{n\rightarrow\infty}\frac{1}{n}S_n\tilde{\zeta}(\om)=\int \tilde{\zeta} d\tilde{\mu}_q 
\]
and
\[ 
\lim_{n\rightarrow\infty}\frac{1}{n}S_n \tilde{f}(\om)=\int \tilde{f} d\tilde{\mu}_q. 
\] 
Therefore, for all $\om\in E_q$ we get 
\[ 
\lim_{n\rightarrow\infty}\frac{S_n \tilde{f}(\om)}{S_n\tilde{\zeta}(\om)}
=\frac{\int \tilde{f} d\tilde{\mu}_q}{\int \tilde{\zeta} d\tilde{\mu}_q}=\a(q). 
\]
Hence $E_q\sbt E_r^\infty(\a(q))\cap E_{rr}^\infty$ and thus
$\tilde{\mu}_q(E_r^\infty(\a(q))\cap E_{rr}^\infty)\geq\tilde{\mu}_q(E_q)=1$.

\

(c) Let $q\in\R$. Using part (a) with $\a=\a(q)$ and part (b), we deduce that 
\[ \mu_q\bigl(J_{r,m_F}(\a(q))\bigr)
\geq\tilde{\mu}_q\circ\pi^{-1}\bigl(\pi(E_r^\infty(\a(q))\cap E_{rr}^\infty)\bigr)
\geq\tilde{\mu}_q\bigl(E_r^\infty(\a(q))\cap E_{rr}^\infty\bigr)
=1.\]
\endpf

\

Let us now remind the reader about the Legendre transform. 
Let $k$ be a strictly convex function on an interval $I$
(hence $k''>0$ wherever this second derivative exists). 
The Legendre transform of $k$ is the function $l$ defined by 
$l(p)=\max\{p\,x-k(x)\}$ wherever the maximum exists. It can 
be proved that the domain of $l$ is either a point, an interval 
or a half-line. It can further be shown that $l$ is strictly 
convex and that the Legendre transform is involutive. We then 
say that the functions $k$ and $l$ form a Legendre transform 
pair. The following theorem (see~\cite{rock}) gives a useful 
characterization of a Legendre transform pair. 

\bthm\label{rockyou}
Two strictly convex differentiable functions $k$ and $l$ form a Legendre transform pair if and only if 
$l(-k'(q))=k(q)-qk'(q)$.
\ethm

We shall now prove that $f_{r,m_F}(\a)$ and $T(q)$ form a Legendre transform pair.
Recall that $T'(q)=-\a(q)$ by Theorem~\ref{alphat}.

\bthm\label{1mam17}
For every $q\in\R$ we have $f_{r,m_F}(\a(q))=q\a(q)+T(q)$. In other terms, $f_{r,m_F}(-T'(q))=T(q)-qT'(q)$.
\ethm

Proof. Using Theorem~\ref{p2mam11}(a,b), 
Theorem~4.4.2 in~\cite{gdms}, Theorem~2.2.9 in~\cite{gdms} which 
guarantees that $\tilde{\mu}_q$ is an (in fact, the unique) equilibrium 
state for $\tilde{f}_q$, and the fact that $P(q,T(q))=P(F_q)=P(\tilde{f}_q)=0$ 
by definition of the temperature function $T(q)$ in Corollary~\ref{c1mam5}, we obtain  
\begin{eqnarray*}
f_{r,m_F}(\a(q))
&=&\HD(J_{r,m_F}(\a(q)))
\geq\HD(\pi(E_r^\infty(\a(q))\cap E_{rr}^\infty))
\geq\HD(\tilde{\mu}_q\circ\pi^{-1}) \\
&=&\frac{h_{\tilde{\mu}_q}(\sg)}{\chi_{\tilde{\mu}_q}(\sg)}
=\frac{P(\tilde{f}_q)-\int \tilde{f}_q\,d\tilde{\mu}_q}{\chi_{\tilde{\mu}_q}(\sg)}
=\frac{-\int \tilde{f}_q\, d\tilde{\mu}_q}{-\int \tilde{\zeta}\, d\tilde{\mu}_q}
=\frac{\int (q\tilde{f}+T(q)\tilde{\zeta})\, d\tilde{\mu}_q}{\int \tilde{\zeta}\, d\tilde{\mu}_q} \\
&=&q\ \frac{\int \tilde{f}\, d\tilde{\mu}_q}{\int \tilde{\zeta}\, d\tilde{\mu}_q}+T(q)
=q\a(q)+T(q).
\end{eqnarray*}

To prove the other inequality, fix $x\in J_{r,m_F}(\a(q))$. Then there is $\om\in E_r^\infty$ such 
that $\pi(\om)=x$. Let  $0<R<\min\{K^{-1},r_{max}(\om)\}$. Let also $\{n_j\}_{j\in\N}$ be any 
subsequence of the increasing sequence $\{n_j(\om,R)\}$ of all $n$'s such that 
$\pi(\sg^n\om)\in X_{t(\om_n)}\backslash B(\partial X_{t(\om_n)},R)$. 
For every $n\in\N$ we have
\[ 
\varphi_{\om|_n}\bigl(B(\pi(\sg^n\om),R)\bigr)
\sbt B\bigl(\pi(\om),KR|\varphi_{\om|_n}'(\pi(\sg^n\om))|\bigr). 
\] 
Like~(\ref{1mam15}), the conformality of $m_q$ and the bounded variation principle 
for $\tilde{f}_q$ give
\begin{equation}\label{1mam18}
\aligned
m_q\bigl(B(x,KR|\varphi_{\om|_n}'(\pi(\sg^n\om))|)\bigr)
&\geq m_q\Bigl(\varphi_{\om|_n}\bigl(B(\pi(\sg^n\om),R)\cap J_{t(\om_n)}\bigr)\Bigr) \\
&\geq \exp\bigl(\inf_{\tau\in[\om|_n]} S_n\tilde{f}_q(\tau)\bigr) 
      \ m_q\bigl(B(\pi(\sg^n\om),R)\bigr) \\
&\geq M_{m_q}(R) B(\tilde{f}_q) \exp\bigl(S_n \tilde{f}_q(\om)\bigr) \\
\endaligned
\end{equation}
where $B(\tilde{f}_q)$ is a constant of bounded variation (see Lemma~2.3.1 in~\cite{gdms}).
Hence
\begin{equation}\label{1mam19}
\aligned
\frac{\log m_q\bigl(B(x,KR|\varphi_{\om|_n}'(\pi(\sg^n\om))|)\bigr)}
     {\log(KR|\varphi_{\om|_n}'(\pi(\sg^n\om))|)} 
&\leq\frac{\log(M_{m_q}(R) B(\tilde{f}_q))+S_n\tilde{f}_q(\om)}
          {\log(KR) +\log|\varphi_{\om|_n}'(\pi(\sg^n\om))|} \\
&=\frac{\log(M_{m_q}(R) B(\tilde{f}_q))
+q S_n\tilde{f}(\om)+T(q)S_n\tilde{\zeta}(\om)}{\log(KR)+S_n\tilde{\zeta}(\om)} \\
&=\frac{\frac{\log(M_{m_q}(R) B(\tilde{f}_q))}{S_n\tilde{\zeta}(\om)}
        +q\frac{S_n\tilde{f}(\om)}{S_n\tilde{\zeta}(\om)}+T(q)}
       {\frac{\log(KR)}{S_n\tilde{\zeta}(\om)}+1}. \\
\endaligned
\end{equation}

\noindent (So far the estimates are valid for all $n\in\N$.) Like~(\ref{1mam14}), every $y\in J_{t(\om_{n_j})}$ admits a $\tau\in [\om|_{n_j}]$ 
such that $y=\pi(\sg^{n_j}\tau)$, and for such $y$ and $\tau$ we have 
$S_{\om|_{n_j}}F(y)=S_{n_j}\tilde{f}(\tau)$. The conformality of $m_F$, 
the fact that $B(\pi(\sg^{n_j}\om),R)\sbt\Int(X_{t(\om_{n_j})})$, 
and the OSC then give
\begin{eqnarray*}
m_F\bigl(B(x,K^{-1}R|\varphi_{\om|_{n_j}}'(\pi(\sg^{n_j}\om))|)\bigr)
&\leq& m_F\bigl(\varphi_{\om|_{n_j}}(B(\pi(\sg^{n_j}\om),R))\bigr) \\
&\leq& m_F\bigl(\varphi_{\om|_{n_j}}(\Int(X_{t(\om_{n_j})}))\bigr) \\
&=   & m_F\bigl(\varphi_{\om|_{n_j}}(\Int(X_{t(\om_{n_j})}))\cap J_{i(\om)}\bigr) \\
&=   & m_F\Bigl(\varphi_{\om|_{n_j}}\bigl(\Int(X_{t(\om_{n_j})})\cap J_{t(\om_{n_j})}\bigr)\Bigr) \\
&\leq& \exp\bigl(\sup_{y\in \Int(X_{t(\om_{n_j})})\cap J_{t(\om_{n_j})}} S_{\om|_{n_j}}F(y)\bigr) 
       \ m_F\bigl(\Int(X_{t(\om_{n_j})})\cap J_{t(\om_{n_j})}\bigr) \\
&\leq&\exp\bigl(\sup_{\tau\in[\om|_{n_j}]} S_{n_j}\tilde{f}(\tau)\bigr) \\
&\leq& B(\tilde{f}) \exp\bigl(S_{n_j}\tilde{f}(\om)\bigr),
\end{eqnarray*}
where $B(\tilde{f})$ is a constant of bounded variation for $\tilde{f}$ (see the bounded 
variation principle for H\"older continuous potentials, Lemma~2.3.1 in~\cite{gdms}). Thus,
\begin{equation}\label{1mam21}
\frac{\log m_F\bigl(B(x,K^{-1}R|\varphi_{\om|_{n_j}}'(\pi(\sg^{n_j}\om))|)\bigr)}
     {\log(K^{-1}R|\varphi_{\om|_{n_j}}'(\pi(\sg^{n_j}\om))|)}
\geq
\frac{\log(B(\tilde{f}))+S_{n_j}\tilde{f}(\om)}{\log(K^{-1}R)+S_{n_j}\tilde{\zeta}(\om)} 
\end{equation}
Using~(\ref{1mam19}) and~(\ref{1mam21}), we deduce that 
\begin{eqnarray*}
d_{m_q}(x)=\lim_{r\rightarrow 0} \frac{\log m_q(B(x,r))}{\log r}
&\leq&\limsup_{j\rightarrow\infty}
\frac{\log m_q\bigl(B(x,KR|\varphi_{\om|_{n_j}}'(\pi(\sg^{n_j}\om))|)\bigr)}
     {\log(KR|\varphi_{\om|_{n_j}}'(\pi(\sg^{n_j}\om))|)} \\
&\leq&\limsup_{j\rightarrow\infty}
\Bigl(q\frac{S_{n_j}\tilde{f}(\om)}{S_{n_j}\tilde{\zeta}(\om)}+T(q)\Bigr) \\
&=&q\limsup_{j\rightarrow\infty}\frac{\log(B(\tilde{f}))+S_{n_j}\tilde{f}(\om)}
                                     {\log(K^{-1}R)+S_{n_j}\tilde{\zeta}(\om)}
                                +T(q) \\
&\leq&q\limsup_{j\rightarrow\infty}
\frac{\log m_F\bigl(B(x,K^{-1}R|\varphi_{\om|_{n_j}}'(\pi(\sg^{n_j}\om))|)\bigr)}
     {\log(K^{-1}R|\varphi_{\om|_{n_j}}'(\pi(\sg^{n_j}\om))|)}
+T(q) \\
&=&q\a(q)+T(q).
\end{eqnarray*}

\noindent As $d_{m_q}(x)\leq q\a(q)+T(q)$ for every $x\in J_{r,m_F}(\a(q))$, 
we deduce that $f_{r,m_F}(\a(q))\leq q\a(q)+T(q)$. 
\endpf

\

All of the above results give us an analog of Theorem~4.9.4 in~\cite{gdms}.
 
\bthm\label{multanalr} 
Let $S=\{\varphi_e:X_{t(e)}\rightarrow X_{i(e)}\,|\,e\in E\}$ be a 
cofinitely regular CGDMS satisfying SOSC and having an underlying 
finitely primitive edge incidence matrix $A$. 
Let $\th$ be the finiteness parameter of $S$.
Suppose that $h_{\tilde{\mu}_q}(\sg)/\chi_{\tilde{\mu}}(\sg)>\th$. Then the following statements hold.
\begin{itemize}
\item[(a)] The number $d_{\mu_F}(x)$ exists for $\mu_F$-a.e. $x\in J_r$ and 
\[ 
d_{\mu_F}(x)
=\frac{\int \tilde{f} \, d_{\tilde{\mu}_{\tilde{f}}}}{\int \tilde{\zeta}\, d_{\tilde{\mu}_{\tilde{f}}}}. 
\]
\item[(b)] The function $T:\R\rightarrow\R$ is real-analytic, $T(0)=\HD(J)$, 
           and $T'(q)<0$, $T''(q)\geq 0$ for all $q\in\R$. 
\item[(c)] For every $q\in\R$, we have $f_{r,\mu_F}(-T'(q))=f_{r,\mu_F}(\a(q))=q\a(q)+T(q)=T(q)-qT'(q)$.  
           That is, $f_{r,\mu_F}(\a)$ and $T(q)$ form a Legendre pair of functions. 
\item[(d)] If $\tilde{\mu}_{\tilde{f}}\neq\tilde{\mu}_{\HD(J)\tilde{\zeta}}$ or, 
           equivalently, if $\tilde{f}$ and $\HD(J)\tilde{\zeta}$ are not 
           cohomologous modulo any constant, then the function
           $\a\mapsto f_{r,\mu_F}(\a)$, $\a\in(\a_1,\a_2)$ is real-analytic, where the interval $(\a_1,\a_2)$, 
           $0\leq\a_1\leq\a_2\leq \infty$, is the range of $-T'(q)$. Otherwise, $T'(q)=\HD(J)$ for every 
           $q\in\R$. 
\end{itemize}
\ethm

{\sl Proof.} (a) Using Theorem~\ref{alphat}(b), notice that $\mu_1=\mu_F$. 
Thus, by Theorem~\ref{p2mam11}(c), we have $\mu_F(J_{r,m_F}(\a(1)))=1$. 
Since $\a(1)=\int \tilde{f}\, d\tilde{\mu}_{\tilde{f}}/\int \tilde{\zeta}\, 
d\tilde{\mu}_{\tilde{f}}$ and since 
the measures $\mu_F$ and $m_F$ are boundedly equivalent, we hence have that 
\[ 
\mu_F\Bigl(J_{r,\mu_F}\bigl(\int \tilde{f}\, d\tilde{\mu}_{\tilde{f}}/
\int \tilde{\zeta}\, d\tilde{\mu}_{\tilde{f}}\bigr)\Bigr)=1. 
\] 
Part (b) is essentially Theorem~\ref{alphat}. 
Part (c) corresponds to Theorem~\ref{alphat} and Theorem~\ref{1mam17}.
Finally, part (d) is a consequence of Lemma~4.9.5 in~\cite{gdms} (with $\Delta_1=\R$) 
and parts (c) and (b) of the present theorem.
\endpf

\section{Multifractal analysis over $J$}\label{jcond}

\subsection{Multifractal analysis over $J$ under the Boundary Separation Condition}\label{ssc}

If $S$ satisfies the Boundary Separation Condition (BSC), that is, if 
\[ \dist(\partial X,\overline{\cup_{i\in I}\varphi_i(X)})>0, \]
then $E_r^\infty=E_A^\infty$
and thus $J_r=J$. Thus, section~\ref{jr} gives us the multifractal analysis over $J$. Indeed,
denoting the set of points of $J$ at which the local dimension of a measure $\mu$ is equal to $\a$ by
$J_\mu(\a)$ and the Hausforff dimension of $J_\mu(\a)$ by $f_\mu(\a)$, Theorem~\ref{multanalr} reduces to the following.  

\bthm\label{multanal} 
Let $S=\{\varphi_e:X_{t(e)}\rightarrow X_{i(e)}\,|\,e\in E\}$ be a cofinitely regular 
CGDMS satisfying SOSC and BSC, and having an underlying finitely 
primitive edge incidence matrix $A$. 
Let $\th$ be the finiteness parameter of $S$.
Suppose that $h_{\tilde{\mu}_q}(\sg)/\chi_{\tilde{\mu}}(\sg)>\th$. Then the following statements hold.
\begin{itemize}
\item[(a)] The number $d_{\mu_F}(x)$ exists for $\mu_F$-a.e. $x\in J$ and 
\[ 
d_{\mu_F}(x)
=\frac{\int \tilde{f} \, d_{\tilde{\mu}_{\tilde{f}}}}{\int \tilde{\zeta}\, d_{\tilde{\mu}_{\tilde{f}}}}. 
\]
\item[(b)] The function $T:\R\rightarrow\R$ is real-analytic, $T(0)=\HD(J)$, 
           and $T'(q)<0$, $T''(q)\geq 0$ for all $q\in\R$. 
\item[(c)] For every $q\in\R$, we have $f_{\mu_F}(-T'(q))=f_{\mu_F}(\a(q))=q\a(q)+T(q)=T(q)-qT'(q)$.  
           That is, $f_{\mu_F}(\a)$ and $T(q)$ form a Legendre pair of functions. 
\item[(d)] If $\tilde{\mu}_{\tilde{f}}\neq\tilde{\mu}_{\HD(J)\tilde{\zeta}}$ or, 
           equivalently, if $\tilde{f}$ and $\HD(J)\tilde{\zeta}$ are not 
           cohomologous modulo any constant, then the function
           $\a\mapsto f_{\mu_F}(\a)$, $\a\in(\a_1,\a_2)$ is real-analytic, where the interval $(\a_1,\a_2)$, 
           $0\leq\a_1\leq\a_2\leq \infty$, is the range of $-T'(q)$. Otherwise, $T'(q)=\HD(J)$ for every 
           $q\in\R$. 
\end{itemize}
\ethm  

An important family of CGDMSs (in fact, conformal IFSs (CIFSs)) which satisfies the boundary 
separation condition are real continued fractions with the digit $1$ deleted.

\bex
Let $-1/4\leq\e<0$. Set $X=[-\e,3/4]$. Let $S=\{\varphi_n:X\rightarrow X\,|\, n\in\N\backslash\{1\}\}$,
where
\[ 
\varphi_n(x)=\frac{1}{n+x}.
\]
Then $S$ is a cofinitely regular CIFS which satisfies 
both the strong open set condition and the boundary separation condition.
\eex

\subsection{Multifractal analysis over $J$ under other conditions}\label{oc}


%
%
%


We shall now prove that $f_{m_F}(\a)$ and $T(q)$ form a Legendre transform pair
under some conditions. Recall that $T'(q)=-\a(q)$ by Theorem~\ref{alphat}.

\bthm\label{cond}
Suppose there exists a countable set $J_0\sbt J$ such that for every 
$x\in J\backslash J_0$ there are $\om\in E_A^\infty$ with $\pi(\om)=x$, a constant 
$C=C(x,\om)>0$, an increasing subsequence 
$\{n_j\}_{j\in\N}=\{n_j(x,\om)\}_{j\in\N}$ of natural numbers and a sequence 
$\{r_{n_j}\}_{j\in\N}=\{r_{n_j}(x,\om)\}_{j\in\N}$ of positive real numbers such that
\begin{enumerate}
\item $m_F(B(x,r_{n_j}))\leq C\exp(S_{n_j}\tilde{f}(\om))$ for all $j\in\N$;
\item $1>r_{n_j}\geq C\|\varphi_{\om|_{n_j}}'\|$ for all $j\in\N$;
\item $\lim_{j\rightarrow\infty}r_{n_j}=0$.
\end{enumerate}
Then $f_{m_F}(\a(q))=q\a(q)+T(q)$ for 
all $q\in\R$.
\ethm

{\sl Proof.} Clearly, $f_{m_F}(\a(q))\geq f_{r,m_F}(\a(q))=q\a(q)+T(q)$ using Theorem~\ref{1mam17}. 
We shall now prove the other inequality. Since $J_0$ is countable, it is sufficient 
to show that $d_{m_q}(x)\leq q\a(q)+T(q)$ for every $x\in J_{m_{F}}(\a(q))\backslash J_0$. 
Accordingly, fix $x\in J_{m_{F}}(\a(q))\backslash J_0$. Let $\om$, $C$,  
$\{n_j\}_{j\in\N}$ and $\{r_{n_j}\}_{j\in\N}$ be as above. 
For every $n\in\N$ we have
\[ 
\varphi_{\om|_n}\bigl(B(\pi(\sg^n\om),R)\bigr)
\sbt B\bigl(\pi(\om),KR|\varphi_{\om|_n}'(\pi(\sg^n\om))|\bigr). 
\] 

\noindent Like in~(\ref{1mam18}), the conformality of $m_q$ and the bounded variation 
principle for $\tilde{f}_q$ give

\[
m_q\bigl(B(x,KR|\varphi_{\om|_n}'(\pi(\sg^n\om))|)\bigr)
\geq M_{m_q}(R) B(\tilde{f}_q) \exp\bigl(S_n \tilde{f}_q(\om)\bigr) \\
\]
where $B(\tilde{f}_q)$ is a constant of bounded variation (see Lemma~2.3.1 in~\cite{gdms}).
Hence, as in~(\ref{1mam19}),

\begin{equation}\label{1mam23}
\aligned
\frac{\log m_q\bigl(B(x,KR|\varphi_{\om|_n}'(\pi(\sg^n\om))|)\bigr)}
     {\log(KR|\varphi_{\om|_n}'(\pi(\sg^n\om))|)} 
&\leq \frac{\frac{\log(M_{m_q}(R) B(\tilde{f}_q))}{S_n\tilde{\zeta}(\om)}
            +q\frac{S_n\tilde{f}(\om)}{S_n\tilde{\zeta}(\om)}+T(q)}
           {\frac{\log(KR)}{S_n\tilde{\zeta}(\om)}+1}. \\
\endaligned
\end{equation}

\noindent So far the estimates are valid for all $n\in\N$. Now, by assumptions~(1) and~(2) 
we obtain for all $j\in\N$
\begin{equation}\label{1mam24}
\aligned
\frac{\log m_F(B(x,r_{n_j}))}{\log r_{n_j}}
&\geq\frac{\log(C)+S_{n_j}\tilde{f}(\om)}{\log r_{n_j}} \\
&\geq\frac{\log(C)+S_{n_j}\tilde{f}(\om)}{\log(C)+\log\|\varphi_{\om|_{n_j}}'\|} \\
&\geq\frac{\log(C)+S_{n_j}\tilde{f}(\om)}{\log(C)+\log|\varphi_{\om|_{n_j}}'(\om)|} \\
&=\frac{\log(C)+S_{n_j}\tilde{f}(\om)}{\log(C)+S_{n_j}\tilde{\zeta}(\om)}.
\endaligned
\end{equation}
Using~(\ref{1mam23}) and~(\ref{1mam24}), we deduce that 
\begin{eqnarray*}
\underline{d}_{m_q}(x)=\liminf_{r\rightarrow 0} \frac{\log m_q(B(x,r))}{\log r}
&\leq&\limsup_{j\rightarrow\infty}\frac{\log m_q\bigl(B(x,KR|\varphi_{\om|_{n_j}}'(\pi(\sg^{n_j}\om))|)\bigr)}
                                       {\log(KR|\varphi_{\om|_{n_j}}'(\pi(\sg^{n_j}\om))|)} \\
&\leq&\limsup_{j\rightarrow\infty}\Bigl(q\frac{S_{n_j}\tilde{f}(\om)}{S_{n_j}\tilde{\zeta}(\om)}+T(q)\Bigr) \\
&=&q\limsup_{j\rightarrow\infty}\frac{\log(C)+S_{n_j}\tilde{f}(\om)}{\log(C)+S_{n_j}\tilde{\zeta}(\om)}+T(q) \\
&\leq&q\limsup_{j\rightarrow\infty}\frac{\log m_F(B(x,r_{n_j}))}{\log r_{n_j}}+T(q) \\
&=&q\a(q)+T(q).
\end{eqnarray*}

\noindent As $\underline{d}_{m_q}(x)\leq q\a(q)+T(q)$ for every $x\in J_{m_F}(\a(q))\backslash J_0$ and $J_0$ is countable, we deduce that 
$f_{m_F}(\a(q))\leq q\a(q)+T(q)$. 
\endpf

\

In the framework of Theorem~\ref{cond}, Theorem~\ref{multanalr} reduces to the following.  

\bthm\label{multanalcond} 
Let $S=\{\varphi_e:X_{t(e)}\rightarrow X_{i(e)}\,|\,e\in E\}$ be a cofinitely regular 
CGDMS satisfying SOSC and BSC, and having an underlying finitely 
primitive edge incidence matrix $A$. 
Let $\th$ be the finiteness parameter of $S$.
Suppose that the conditions of Theorem~\ref{cond} are fulfilled. Suppose also that $h_{\tilde{\mu}_q}(\sg)/\chi_{\tilde{\mu}}(\sg)>\th$. Then the following statements hold.
\begin{itemize}
\item[(a)] The number $d_{\mu_F}(x)$ exists for $\mu_F$-a.e. $x\in J_r$ and 
\[ 
d_{\mu_F}(x)
=\frac{\int \tilde{f} \, d_{\tilde{\mu}_{\tilde{f}}}}{\int \tilde{\zeta}\, d_{\tilde{\mu}_{\tilde{f}}}}. 
\]
\item[(b)] The function $T:\R\rightarrow\R$ is real-analytic, $T(0)=\HD(J)$, 
           and $T'(q)<0$, $T''(q)\geq 0$ for all $q\in\R$. 
\item[(c)] For every $q\in\R$, we have $f_{\mu_F}(-T'(q))=f_{\mu_F}(\a(q))=q\a(q)+T(q)=T(q)-qT'(q)$.  
           That is, $f_{\mu_F}(\a)$ and $T(q)$ form a Legendre pair of functions. 
\item[(d)] If $\tilde{\mu}_{\tilde{f}}\neq\tilde{\mu}_{\HD(J)\tilde{\zeta}}$ or, 
           equivalently, if $\tilde{f}$ and $\HD(J)\tilde{\zeta}$ are not 
           cohomologous modulo any constant, then the function
           $\a\mapsto f_{\mu_F}(\a)$, $\a\in(\a_1,\a_2)$ is real-analytic, where the interval $(\a_1,\a_2)$, 
           $0\leq\a_1\leq\a_2\leq \infty$, is the range of $-T'(q)$. Otherwise, $T'(q)=\HD(J)$ for every 
           $q\in\R$. 
\end{itemize}
\ethm  

\

There are interesting families of CGDMSs (in fact, even of CIFSs) which satisfy the conditions 
imposed in Theorem~\ref{cond}. Among others, let us mention real continued fractions over the interval 
$[0,1]$.

\bex
Let $X=[0,1]$. Let $S=\{\varphi_n:X\rightarrow X\,|\, n\in\N\}$,
where
\[ 
\varphi_n(x)=\frac{1}{n+x}.
\]
Then $S$ is a cofinitely regular CIFS which satisfies 
the SOSC and conditions (1)--(3) of Theorem~\ref{cond} 
(but not the boundary separation condition).
Moreover, writing 
\[ x=\frac{1}{x_1+\frac{1}{x_2+\ldots}}\]
we note that 
$J_r=\{x\in[0,1]:\liminf_{n\rightarrow\infty}x_n<\infty\}$ and thus $J_r^c=\{x\in[0,1]:\lim_{n\rightarrow\infty}x_n=\infty\}$.
\eex

In fact, there is a larger class of one-dimensional CIFSs for which 
conditions (1)--(3) of Theorem~\ref{cond} are fulfilled. By one-dimensional, 
we simply mean that $X$ is a subinterval of $\R$.

\bthm
Let $S=\{\varphi_n:X\rightarrow X\,|\,n\in\N\}$ be a one-dimensional CIFS 
satisfying the following properties. 
\begin{itemize}
\item[(i)] $\|\varphi_n'\|$ is comparable to $\|\varphi_{n+1}'\|$, that is, 
there exists a constant $C\geq 1$ such that 
\[
C^{-1}\leq\frac{\|\varphi_{n+1}'\|}{\|\varphi_n'\|}\leq C, \hspace{1cm} \forall\,\,n\in\N.
\]
\item[(ii)] The generators $\varphi_n$, $n\in\N$, of $S$ either all preserve orientation 
(i.e. $\varphi_n'>0$ for all $n\in\N$) or all reverse orientation (i.e. $\varphi_n'<0$ for all $n\in\N$).
\item[(iii)] Either $\varphi_{n+1}>\varphi_n$ for all $n\in\N$ or $\varphi_{n+1}<\varphi_n$ for all $n\in\N$.
\end{itemize}
Then all three conditions in Theorem~\ref{cond} are met.
\ethm 

{\sl Proof.} It follows from the bounded distortion property and (i) that 
\begin{equation}\label{compdiam}
(CK)^{-1}\leq\frac{|\varphi_{n+1}(X)|}{|\varphi_n(X)|}\leq CK, \hspace{1cm} \forall\,\,n\in\N.
\end{equation}
Moreover, a rather straightforward calculation shows that the amalgamated function 
$\tilde{\zeta}$ satisfies 
\[ 
\sup_{\om,\tau\in E_A^\infty:|\om_1-\tau_1|=1}|\tilde{\zeta}(\om)-\tilde{\zeta}(\tau)|\leq\log(CK):=D. 
\]  
%
Now, let 
\[ J_0=\bigl\{\pi(\tau)\,\bigl|\,\sg^k\tau=1^\infty \mbox{ for some } k\in\N \bigr\}. \]
Clearly, $J_0$ is countable. Let $x\in J\backslash J_0$ and $\om\in\N^\infty$ 
such that $\pi(\om)=x$. Let $n_j$ be the $j$-th letter in the word $\om$ which 
is not a $1$. Thus, $\om_{n_j}\neq 1$. Let $r_{n_j}=|\varphi_{\om_{n_j}}(X)|/(CK)$.
Therefore $r_{n_j}\leq\min\{|\varphi_{\om_{n_j}-1}(X)|,|\varphi_{\om_{n_j}}(X)|,|\varphi_{\om_{n_j}+1}(X)|\}$ 
by~(\ref{compdiam}) and hence 
\[ B(x,r_{n_j})\cap J\sbt\bigcup_{k=\om_{n_j}-1}^{\om_{n_j}+1}\varphi_{\om|_{n_{j}-1}k}(X)\cap J. \] 
Then, since $P(F)=0$, 
\begin{eqnarray*}
m_{F}\bigl(B(x,r_{n_j})\bigr)
&\leq&\sum_{k=\om_{n_j}-1}^{\om_{n_j}+1}m_{F}\bigl(\varphi_{\om|_{n_j-1}k}(X)\bigr) \\
&\leq&\sum_{k=\om_{n_j}-1}^{\om_{n_j}+1}
      \int_{J}\exp\bigl(S_{\om|_{n_j-1}k}F(y)-n_jP(F)\bigr)dm_{F} \\
&\leq&\sum_{k=\om_{n_j}-1}^{\om_{n_j}+1}
      \exp\bigl(\sup_{y\in J}S_{\om|_{n_j-1}k}F(y)\bigr) \\
&=&\sum_{k=\om_{n_j}-1}^{\om_{n_j}+1}
      \exp\bigl(\sup_{\tau\in[\om|_{{n_j}-1}k]}S_{n_j}\tilde{f}(\tau)\bigr) \\
&\leq&B(\tilde{f})\sum_{k=\om_{n_j}-1}^{\om_{n_j}+1}
      \exp\bigl(S_{n_j}\tilde{f}(\om|_{{n_j}-1}k\om|_{n_j+1}^\infty)\bigr) \\
&=&B(\tilde{f})\sum_{k=\om_{n_j}-1}^{\om_{n_j}+1}
   \exp\bigl(S_{n_j-1}\tilde{f}(\om|_{{n_j}-1}k\om|_{n_j+1}^\infty)+\tilde{f}(k\om|_{n_j+1}^\infty)\bigr) \\
&=&B(\tilde{f})\sum_{k=\om_{n_j}-1}^{\om_{n_j}+1}
   \exp\bigl(S_{n_j-1}\tilde{f}(\om|_{{n_j}-1}k\om|_{n_j+1}^\infty)-S_{n_j-1}\tilde{f}(\om)\bigr)
   \exp\bigl(S_{n_j-1}\tilde{f}(\om)\bigr) \\
&&\hspace{5cm}\cdot\exp\bigl(\tilde{f}(k\om|_{n_j+1}^\infty)-\tilde{f}(\om|_{n_j}^\infty)\bigr)
                   \exp\bigl(\tilde{f}(\om|_{n_j}^\infty)\bigr) \\
&\leq&B(\tilde{f})\sum_{k=\om_{n_j}-1}^{\om_{n_j}+1}B(\tilde{f})
      \exp\bigl(S_{n_j-1}\tilde{f}(\om)\bigr)\exp(|k-\om_{n_j}|D)
      \exp(\tilde{f}(\om|_{n_j}^\infty)) \\
&\leq&3e^D(B(\tilde{f}))^2\exp\bigl(S_{n_j}\tilde{f}(\om)\bigr).
\end{eqnarray*}
Thus, condition (1) in Theorem~\ref{cond} is satisfied. Note also that condition (2) is fulfilled since 
\[ 
r_{n_j}=\frac{1}{CK}|\varphi_{\om_{n_j}}(X)|\geq\frac{1}{CK^2}\|\varphi_{\om_{n_j}}'\|. 
\]
Finally, condition (3) is obviously satisfied, as $\sum_{j\in\N}|\varphi_{\om_{n_j}}(X)|\leq|X|<\infty$ and hence $r_{n_j}\rightarrow 0$.
\endpf 





\end{document}